\documentclass[onefignum,onetabnum]{siamonline190516}

\usepackage[firstpage]{draftwatermark}
\SetWatermarkText{PREPRINT - To be published in SIAM Journal on Applied Dynamical Systems.}
\SetWatermarkFontSize{0.4cm}
\SetWatermarkScale{1}
\SetWatermarkAngle{0}
\SetWatermarkColor[gray]{0.4}
\SetWatermarkVerCenter{2.5cm}

\usepackage{lipsum}
\usepackage{amsfonts}
\usepackage{graphicx}
\usepackage{epstopdf}
\usepackage{algorithmic}
\ifpdf
  \DeclareGraphicsExtensions{.eps,.pdf,.png,.jpg}
\else
  \DeclareGraphicsExtensions{.eps}
\fi

\newsiamremark{remark}{Remark}
\newsiamremark{hypothesis}{Hypothesis}
\crefname{hypothesis}{Hypothesis}{Hypotheses}
\newsiamthm{claim}{Claim}

\usepackage{mathrsfs} 
\usepackage[RPvoltages]{circuitikz}
\usepackage{enumitem}
\setlist[enumerate]{leftmargin=.5in}
\setlist[itemize]{leftmargin=.5in}  
\usepackage{mathtools} 

\DeclareMathOperator{\im}{im}
\DeclareMathOperator{\In}{In}
\colorlet{My_green10}{green!50!black}

\headers{Nonlinear consensus on networks}{M. Homs-Dones,  K. Devriendt and R.  Lambiotte}

\title{Nonlinear consensus  on networks: equilibria, effective resistance and  trees of motifs\thanks{PREPRINT - To be published in SIAM Journal on Applied Dynamical Systems.
\funding{M.H. was supported by ``la Caixa''
Foundation (ID 100010434) with fellowship code LCF/BQ/ EU20/11810061. K.D. was supported by The Alan Turing Institute under the EPSRC
grant EP/N510129/1.}}}

\author{Marc Homs-Dones\thanks{Mathematics Institute, University of Warwick, UK
  (\email{marc.homs-dones@warwick.ac.uk}).}
\and Karel Devriendt\thanks{Mathematical Institute, University of Oxford, UK 
  (\email{devriendt@maths.ox.ac.uk}).}
  \and Renaud Lambiotte\thanks{Mathematical Institute, University of Oxford, UK 
  (\email{renaud.lambiotte@maths.ox.ac.uk}, \url{https://www.maths.ox.ac.uk/people/renaud.lambiotte}).}}

\usepackage{amsopn}

\ifpdf
\hypersetup{
 pdftitle={Nonlinear consensus  on networks: equilibria, effective resistance and  trees of motifs},
 pdfauthor={Marc Homs-Dones, Karel Devriendt and Renaud Lambiotte}
}
\fi

\begin{document}
  
\maketitle 










 









\begin{abstract}
 We study a generic family of nonlinear dynamics on undirected networks  generalising linear consensus. We find a compact expression for its equilibrium points in terms of the  topology of the network and classify their stability using the effective resistance of the underlying graph equipped with appropriate weights. Our general results are  applied to some specific networks, namely trees, cycles and complete graphs. When a network is formed by the union of two subnetworks joined in a single node, we show that the equilibrium points and stability in the whole network can be found by simply studying the smaller subnetworks instead. Applied recursively, this property opens the possibility to investigate the dynamical behaviour on families of networks made of trees of motifs.
\end{abstract}



\section{Introduction}

A broad range of systems can be represented by dynamical systems on networks \cite{newman2018networks}. In this framework, each node is endowed with a time-varying state whose dynamics depend on its own state and the states of its neighbours. Important examples include models of collective behaviour where the decision process is distributed rather than centralised \cite{blondel2005convergence}. 
 As reviewed in \cite{review_app,4118472}, applications can be found in a variety of disciplines such as ecology, where it is used to model animal behaviour and flocking 
\cite{IainD.Couzin2005Elad,TannerH.G2003Sfom,1605401}.
In engineering, it is essential for the design of decentralized control strategies for the movement of robots \cite{LawtonJ.R.T2003Adat}, the rendezvous problem \cite{1673588} and to coordinate decision making when multiple nosy sensors detect an event \cite{AlanyaliM2004DBht}. In the social sciences, it is used as a first model for opinion and language dynamics \cite{CastellanoClaudio2009Spos}, and can be applied in economics to coordinate a decentralized network of buyers \cite{BausoD2003Dcpf}. 
The simplest algorithms and models for consensus are variations of the DeGroot model \cite{degroot1974reaching}, usually called linear consensus dynamics, where the state of a node evolves towards the average value of its neighbours. These models are theoretically appealing thanks to their simplicity, as their dynamics is entirely determined by the spectral properties of the coupling matrix between the elements. However, linear models often arise as approximations for more complicated coupling functions and are not always sufficient to explain their complex behaviour \cite{IainD.Couzin2005Elad}, which calls for the study of nonlinear consensus \cite{Onbif}.
 In many of these applications, nonlinearity has thus been introduced \cite{Arcak_2007,BausoD2003Dcpf}, the most prominent example of nonlinear process being the Kuramoto model of coupled oscillators \cite{book_kuramoto,ArenasA.2007Sicn}, with applications in physics, chemistry, biology and engineering \cite{DorflerFlorian2014Sicn,NabetBenjamin2009DoDM, Tyloo_resistances_in_Kuramoto}.

The focus of this work is on studying the properties of a general class of nonlinear consensus dynamics taking place on a network. We derive exact results for arbitrary choices of undirected networks and a general class of nonlinear consensus dynamics, in contrast with  most other works that focus on mean-field approximations,  specific types of networks or specific classes of nonlinear dynamics.  In particular, we manage to extend the theoretical analysis of \cite{karel2020arXiv}, where the purely graph-theoretic notion of effective resistance was introduced to study a particular nonlinear system. The role of effective resistances in nonlinear dynamical systems was also highlighted in \cite{Tyloo_resistances_in_Kuramoto}, where resistances were shown to be a key determinant of the robustness of a nonlinear system to external perturbations. Another relevant line of research into generic nonlinear network systems can be found in \cite{MartinGolubitsky2006Ndon,GandhiPunit2020BoFI}, where it is shown that certain properties (symmetries) of the underlying network restrict the possible dynamical features (equilibria, periodic states) of any possible system on the network. 
Finally, the work of Bronski and DeVille \cite{BronskiJaredC.2014STfD} presents complementary stability results based on the  study of spanning trees and homology theories of signed graphs, which coincide with ours in simple cases, as we will show in the manuscript.

Our problem is formalised as follows.
Let $n\in \mathbb N$ and $  \mathcal G=(\mathcal V,\mathcal E)$ be a connected undirected network with $\mathcal V=\{1,\dots,n\}$. If $\{i,j\}\in \mathcal E$ we say they are adjacent and write $i\sim j$. Note that we consider only networks without self-loops, so that $i\not \sim i$, as self-loops would not change the dynamics that we will study. Linear consensus dynamics on $\mathcal G$ is defined over the state space $\mathbf x \in \mathbb R^n$, where each node has a scalar state $x_i$, by
\[\dot x_i=-\sum_{j\sim i}(x_i-x_j) \hspace{20pt}\textrm{for all}\hspace{20pt} i\in \mathcal V. \]
Compactly we may write $\dot{\mathbf x }= -L\mathbf x$ where $L$  is the Laplacian matrix of  $\mathcal G$ (see \cref{eq:ddT=L}).  As presented in  \cite{Onbif}, there are three obvious ways to extend this simple linear model by introducing a function $\hat f:\mathbb R \rightarrow \mathbb R$,
 \begin{equation}
     \dot x_i=\sum_{j\sim i} \left (\hat f(x_i)-\hat f(x_j)\right),\hspace{0.45cm} \dot x_i=\sum_{j\sim i} \hat f\left (x_i-x_j\right),\hspace{0.45cm}\dot x_i=\hat f\left (\sum_{j\sim i} \left (x_i-x_j\right)\right ).
     \label{eq:three_ODE}
 \end{equation}
 When $\hat{f}=-\textrm{id}$, all three cases recover the linear consensus model. If we denote by $f:\mathbb R^n\rightarrow \mathbb R^n  $ the function that applies $\hat f$ to each component, we may write the first and third ODE above more compactly as
 \begin{equation}
 \label{eq:compact:form}
\dot{\mathbf x }= Lf(\mathbf x) \hspace{1cm} \textrm{and}\hspace{1cm}\dot{\mathbf x }=  f(L\mathbf x).  \end{equation}
 As shown in \cite{Onbif} for the case $\hat f(y)=\lambda y-y^3$, this matrix formulation facilitates the theoretical understanding of these systems. In this work we will focus on the second ODE, which presents more challenges as it lacks an obvious matrix formulation. Let $\{f_e\}_{e\in \mathcal E}$ be a family of odd, continuously differentiable\footnote{It would be enough for $f_e$ to be locally Lipschitz as long as they are continuous differentiable in a neighbourhood of the equilibria of \cref{principal}.} functions, where we will denote $f_{\{i,j\}}=f_{ij}=f_{ji}$, and consider the ODE\footnote{Without loss of generality we could always consider the complete graph while choosing  $f_{ij}$ to be null in some edges. We do not choose to do so as in many cases we will require that all $f_{ij}$ coincide. }
 \begin{equation}
\dot x_i=\sum_{j\sim i}f_{ij}(x_i-x_j) \hspace{20pt}\textrm{for all}\hspace{20pt} i\in \mathcal V. 
\label{principal}
\tag{\#}
\end{equation} 
Importantly, this equation allows us to model weighted networks by absorbing the corresponding weights in the functions $f_{ij}$.
As $\mathcal G$ is undirected, we require $f_{ij}$ to be odd so the system exhibits a conservation of quantity, which is desirable in many real scenarios. We may think for instance that a liquid flows between adjacent nodes. Moreover, this choice does not severely limit the applicability of our model, as in many applications a global reference frame is not available \cite{6859235} and thus  we need invariant dynamics  under rotation and translation of  inertial frame. In \cite{7330544} it is shown that this is satisfied if and only if the coupling functions $f_{ij}$ are quasi-linear, i.e. $f_{ij}(y)=k_ {ij}(|y|)y$  for some functions $k_{ij}$, which are a collection of odd functions. 
 It is also worth mentioning that the coupled phase oscillators model, of which the Kuramoto model is a subclass, is defined by
 \[\dot {\theta_i}=w_i-\sum_{j=1}^n a_{\{i,j\}}\sin (\theta_i-\theta_j).\]
 This model is not quite included in \cref{principal}, due to the constants $w_i$. However, the differential of this system coincides with the differential of the case $w_i=0$, which is included in  \cref{principal}. Thus,  most of the presented results regarding linear stability of equilibria  will be applicable to the coupled phase oscillator model. 
 

In this paper, we have obtained four main results for the study of \cref{principal}.
Firstly,  we show that \cref{principal} can be expressed as a matrix product with a nonlinear coupling  \cref{principal_2} which allows us to find a compact expression for the set of equilibrium points  \cref{eq:equlibrium_edges}. This expression relates the equilibria with the topology of the underlying network through the cycle and cut space.
Secondly, we find stability criteria for the equilibrium points using the concept of effective resistance, which in certain situations is tight (see \cref{thm:resistencia:one_edge} and \ref{thm:unstability_r^-_ij r^+_ij}). 
Thirdly, when there are few edges contributing to instability, the Schur complement can be used to reduce the stability problem to a smaller network (see \cref{thm:schur_comp}).  
Similarly, in \cref{section:union_graphs}, we show that if a network consists of two subnetworks joined in a single common node, the task of finding equilibrium points and their stability in the whole network can be reduced to the same task for each of the subnetworks. This result motivates the study of particular classes of networks, i.e. building blocks, that will be combined to form more complex networks,  which we refer to as trees of motifs (see \cref{sec:tree_motifs}). 
In particular, we apply our general results to tree, cycle and complete graphs and derive some additional properties that follow from their specific structures, see \cref{chapter:application_to_network}. For the cycle graph we show in \cref{prop:continuom_equilibria_polinomial} that a particular set of polynomials induce a continuum of equilibria in \cref{principal}. This was observed for a particular case in \cite{Onbif} and \cite{karel2020arXiv} but no explanation for such behaviour was given.

 \section{Mathematical preliminaries} 
 \label{sec:math_prelim}
We give a brief recap of the main structures and results that we will need throughout the paper. 
A \emph{network} or graph $\mathcal G$, which we will always consider undirected,  is given by a pair of sets $(\mathcal V, \mathcal E)$ where $\mathcal V$ is the set of \emph{nodes}, which we will assume of the form $\{1,\dots, n\}$, and $\mathcal E$ a set of unordered, distinct pairs of nodes corresponding to the \emph{edges}. We will sometimes  abuse notation and write $i\in \mathcal G$ or $\{i,j\}\in \mathcal G$ instead of $i\in \mathcal V$ or $\{i,j\}\in \mathcal E$. Given a partition $\mathcal V=\mathcal V_1\cup \mathcal V_2$ of the nodes, we say that the set of edges from $\mathcal V_1$ to $\mathcal V_2$ forms a \emph{cut-set} if it is non-empty.

Associated to a network, we have two  $\mathbb R$ vector spaces  $C_0$, $C_1$ generated by the orthonormal bases $\mathcal V$, $\mathcal E$ respectively. They are related by the \emph{boundary} and coboundary linear maps,
\[d:C_1\longrightarrow C_0, \hspace{1cm} d^\top:C_0\longrightarrow C_1,\]
defined by $d(\{j,i\})=i-j$ where $j<i$. Choosing an ordering in the bases, $d$ is simply the ``signed'' incidence matrix, see   \cref{eq:d_def}. 
The matrices $d$ and $d^\top$ are directly related with consensus dynamics as the \emph{Laplacian} may be defined by (see \cite[Proposition 4.8]{book_Alg_graph}),
\begin{equation}
    L := d\hspace{1pt}d^\top .
    \label{eq:ddT=L}
\end{equation}
Additionally, the kernel and image of these linear maps have a strong topological interpretation \cite{book_Alg_graph}. 
 Indeed, $\ker d^\top$, $\ker d$ and $\im d^\top =(\ker d)^\bot$ are generated by the connected components, the cycles and the cut-sets respectively. In particular, if the network is connected, $\ker d^\top=\langle\mathbf 1\rangle$ and,
 \begin{equation}
     d^\top |_{ \langle \mathbf 1\rangle ^\bot}: \langle \mathbf 1\rangle ^\bot\longrightarrow \im d^\top = (\ker d)^\bot  ,
     \label{eq:isomorphism}
 \end{equation}
is an isomorphism.

A \emph{weighted network} is a network with a set of non-negative  weights $\{w_{\{i,j\}}\}_{\{i,j\}\in \mathcal E}$. Often we only consider the edges with positive weights, for instance  we say that a weighted network is connected, if the set of edges with positive weights is connected in the usual sense.  The Laplacian matrix of a weighted network  with vector of weights $w=(w_{\{i,j\}})_{\{i,j\}\in \mathcal E}$ is defined as 
 \begin{equation}
     L := d \hspace{0.1cm}\textrm{diag}\left (w\right )\hspace{1pt}d^\top,
     \label{eq:laplacian_matrix_prod}
 \end{equation}
which is symmetric and positive semi-definite. Moreover, its null space is generated by the indicator vectors of the connected components of the underlying weighted network. As there is a one to one correspondence between weighted Laplacians and networks, we may directly refer to the connected components of the Laplacian.


Given an ODE, $\dot {\mathbf x}=g(\mathbf x)$, we say that $\mathbf x^\star$ is an \emph{equilibrium point}, if $g(\mathbf x^\star)=0$.
We say that an equilibrium point is \emph{attractive} if solutions that start close enough to it, eventually converge to it.
We say that an equilibrium point is \emph{stable} if solutions that start close enough to it remain close for all positive times, otherwise we say it is \emph{unstable}  (see \cref{supp:sec:dynamics} for formal definitions). In practice it is easier to check if an equilibrium is linearly stable, which implies both being stable and attractive. We say that an equilibrium point is \emph{linearly stable} (resp. \emph{unstable}) if all (resp. any) eigenvalues of the Jacobian matrix of $g$ at that point have negative (resp. positive) real part. Analogously we define linear stability of a matrix. Linear stability also has the advantage  that persists under small perturbations of the vector field \cite{book_Gerald_Teschl}, which is key when the system is an imperfect representation of a phenomena, such as in scientific models.

Given a real symmetric matrix $M$ we define its \emph{inertia}, and denote it by $\In (M)$, as the tuple which entries represent the number of positive, negative and zero eigenvalues (including multiplicities).  This quantity will be useful as it determines the linear stability of the matrix. It is well known that a symmetric matrix can be expressed as $M=\sum_{i=1}^k \lambda_i \mathbf v_i \mathbf v_i^\top$, where $\lambda_i$ are the non-zero eigenvalues and $\mathbf v_i$ are the corresponding eigenvectors forming an orthonormal set. We define the \emph{pseudoinverse} $M^{\dagger}$ of $M$ as
 \[M^\dag :=\sum_{i=1}^k \lambda_i^{-1} \mathbf v_i \mathbf v_i^\top,\]
so $\ker M=\ker M^\dag$ and if  $V\bot \ker M$, then $M^\dag|_V=(M|_V)^{-1}$. It is worth mentioning that we do not need to know the eigenvalues of a matrix to find its pseudoinverse, as it can be computed through column/row operations and matrix mutiplications \cite{full_rank}.



 We denote by $\mathbf 0_k$ (resp. $\mathbf 1_k$) the  column vectors with all entries 0 (resp. 1) and $\mathbf 0_{k,k'}=\mathbf 0_k\mathbf 0_{k'}^\top$. When the dimension is clear we will suppress the subindex. We denote by $\mathbf e_1,\dots, \mathbf e_k$ the standard basis of $\mathbb R^k$ and by $I_k$ the identity matrix. 

 \section{Development of the Model}
 \subsection{First results}
 From now on we study the system \cref{principal} defined on an undirected, connected network $\mathcal G$.
 The following two observations were already mentioned in \cite{karel2020arXiv}. 

The mean state $\overline {\mathbf x }=\frac{1}{n}\sum_{i=1}^n x_i$ is a conserved quantity of the system \cref{principal}. To see this, simply note that the functions $f_{ij}$ are odd, and thus
\[\dot {\overline {\mathbf x }}=\langle \nabla \overline {\mathbf x },\dot{\mathbf x}\rangle=\frac{1}{n}\langle \mathbf 1,\dot{\mathbf x}\rangle=0,\]
i.e. $\overline {\mathbf x }$ is a first integral of \cref{principal}. Thus, the planes defined by $\overline {\mathbf x }=k$ for some constant $k$, which are precisely the affine perpendicular planes to $\mathbf 1$, are invariant. Moreover, as only differences of elements $x_i$ appear in \cref{principal}, the vector field in $\mathbf x$ and $\mathbf x + \lambda \mathbf 1$  is the same for all $\lambda \in \mathbb R$. So we conclude that the dynamics in each plane  are exactly the same and we can simply restrict our study to one of these planes.
From now on, we will study \cref{principal} in the vector space  $\langle \mathbf 1\rangle ^\bot$ or equivalently the invariant plane given by $\overline {\mathbf x } =0$, unless stated otherwise. In \cite{karel2020arXiv} a similar approach was taken, by limiting the study to the state space $\mathbb R^n/\langle\mathbf{1} \rangle$.

The fact that $\overline {\mathbf x }$ is a conserved quantity is also desirable as it shows that our model preserves some of the key properties of linear consensus dynamics. In the linear case $\overline {\mathbf x }$ is not only conserved, but all entries of a solution tend to this value. 
 
 Another important property which is preserved from linear consensus is that the dynamics given by \cref{principal} are gradient dynamics in $\mathbb R^n$.  Indeed, if we let $F_{ij}$ be an antiderivative of $f_{ij}$ for all $\{i,j\}\in \mathcal E$ and define,
 \begin{equation}V(\mathbf x ) :=-\sum_{\{i,j\}\in \mathcal E} F_{ij}(x_i-x_j),
 \label{eq:def_grad}
 \end{equation}
which is well defined as all $F_{ij}$ are even, we get,
 \[\left [-\nabla V(\mathbf x)\right]_i=\sum_{j\sim i}f_{ij}(x_i-x_j)=\dot x_i,\]
 so  $\dot{\mathbf x}= -\nabla V(\mathbf x)$. Thus, $V$ is a potential of the system \cref{principal} and   $V(\mathbf x(t))$ is a decreasing function in time (except in equilibrium points), see \cref{eq:derivative_grad}.  In particular, by LaSalle's invariance principle \cite[Theorem 6.15]{book_Gerald_Teschl} all forward (resp. backwards) bounded orbits ``tend to'' (resp. ``come from'') equilibrium points. In \cref{prop:descomp_state_W} we show that, if $f_{ij}$ satisfy certain conditions, then all forward orbits are bounded and thus converge to equilibria. 

 \subsection{Matrix reformulation}
\label{sec:matrix_reformulation}
Define $f:C_1\longrightarrow C_1$ by 
$f((y_e)_{e\in \mathcal E})=(f_{e}(y_e))_{e\in  \mathcal E}$.
Then, $[ f (d^\top (\mathbf{x}))]_{\{j,k\}}=f_{kj}(x_k-x_j)$ where $\{j,k\}\in \mathcal E$ with $j<k$ and,
\begin{equation}
[d]_{i,\{j,k\}}=\begin{cases}
1		& $ if $ k=i,\\
-1		& $ if $ j=i,\\
 0 &  $ else.$
  \end{cases}
    \label{eq:d_def}
\end{equation}
 So we get, 
\begin{align*}
[d( f(d^\top (\mathbf x)))]_i=&\sum_{\{j,k\}\in \mathcal E }[d]_{i,\{j,k\}} \cdot [ f (d^\top (\mathbf{x}))]_{\{j,k\}}\\=& \sum_{\substack{\{j,i\}\in \mathcal E\\j<i}}f_{ij}(x_i-x_j) + \sum_{\substack{\{i,k\}\in \mathcal E\\i<k}}-f_{ki}(x_k-x_i)=\sum_{j\sim i} f_{ij}(x_i-x_j)=\dot x_i,
\end{align*}
  where in the second to last equality we use that $f_{ki}=f_{ik} $ and that they are odd functions. Hence we can rewrite \cref{principal} as 
 \begin{equation}
    \dot{\mathbf{x}} = (d\circ  f \circ d^\top )(  \mathbf{x}).
     \label{principal_2}
 \end{equation}
This is very reminiscent of the compact forms in equation \cref{eq:compact:form}. Indeed, if we recall \cref{eq:ddT=L}, the three ODEs in \cref{eq:three_ODE} can be expressed as, 
\begin{equation}
    \dot{\mathbf x}=(d\circ d^\top \circ f)(\mathbf x ),\hspace{1cm}\dot{\mathbf x}=(d\circ  f \circ d^\top)(\mathbf x ),\hspace{1cm}\dot{\mathbf x}=(f\circ d\circ d^\top) (\mathbf x ),
\label{eq:compact_form2}
\end{equation}
where in the first and last case $f$ is given by not necessarily odd functions in nodes, i.e. $f((x_v)_{v\in \mathcal V})=(f_{v}(x_v))_{v\in \mathcal V}$ (in \cref{eq:three_ODE} we chose $f_v=\hat f$). In \cite{Onbif} it is shown that the matrix form of the first and third ODE in \cref{eq:compact_form2} is  very useful to study the structure of the corresponding systems. We will show that the same is true for \cref{principal_2}. 

 \Cref{principal_2} also allows us to represent system \cref{principal} as a dynamical system on the edges. That is, in the coordinates  $\mathbf{y}=d^\top \mathbf{x}$,  system \cref{principal} has the form, 
 \begin{equation}
     \dot{\mathbf y}= (d^\top\circ d \circ  f)(\mathbf y),
     \label{principal_y}
 \end{equation}
 and is defined on the cut space $(\ker d)^\bot$, see \cref{eq:isomorphism}. Note that the value in an edge  $\{j,i\}\in\mathcal E$ is given by $x_i-x_j$, where $j<i$.


\subsection{Equilibrium points}
\label{sec:equilibrium_points}
In this section, we find expressions that relate the set of equilibrium points of \cref{principal} with the topology of the underlying network. First, using  \cref{principal_y} and  $\ker (d^\top d)=\ker d$ it is clear that the equilibrium points in the \emph{edge space}, i.e. in the coordinates $\mathbf y$,  are given by,
 \begin{equation}
     \mathscr E_{\mathbf y} := (\ker d)^\bot\cap  f^{-1}(\ker d).
     \label{eq:equlibrium_edges}
 \end{equation}
 Recall that $\ker d$ is the cycle space and $(\ker d)^\bot$ the cut space, which have strong topological interpretations and are easy to compute intuitively. In \cref{supp:sec:dim_equilibria} we further develop this intuition and informally show why we should expect the set of equilibria to be at most $m-n+1$ dimensional, using  \cref{eq:equlibrium_edges}.

 

 The set of equilibrium points is better understood for specific systems such as the Kuramoto model \cite{Dhagash_equilibria_of_Kuramoto, Delabays_equilibria_in_Kuramoto}, but when dealing with complex networks and functions $f_{ij}$, it will generally be unfeasible to find all equilibria explicitly. Thus, in some occasions we will restrict our study to the  detailed-balance stationary states, as  done in \cite{karel2020arXiv}, which are given by,
 \[\tilde {\mathscr E}_{\mathbf y}:= (\ker d) ^\bot \cap f^{-1}(\mathbf 0).\]
In node coordinates, $\mathbf x$, we have, 
\[\mathscr E_{\mathbf x}:=\langle \mathbf 1\rangle ^\bot \cap( d^\top )^{-1}(f^{-1}(\ker d)) ,\]
and,
\[\tilde {\mathscr E}_{\mathbf x}:=\{\mathbf x\in \langle \mathbf 1 \rangle^\bot \textrm{ : } f_{ij}(x_i-x_j)=0 \textrm{ for all } \{i,j\}\in \mathcal E \}.\]
In particular, as odd functions have $0$ as a root, $\mathbf 0$ is an equilibrium point of \cref{principal} in $\langle \mathbf 1\rangle ^\bot$ and thus, $\langle \mathbf 1 \rangle$ is a line of equilibria in $\mathbb R^n$. This is a desired property for generalizations of linear consensus, as in the linear case the equilibrium points are given by the \emph{consensus states}, i.e. $\langle \mathbf 1 \rangle$ or equivalently $x_i=x_j$ for all $i,j$. Note that in our general context we may have other equilibria besides the consensus states. In general this is desirable as it allows for more flexibility, see for instance the animal group decision-making model with non-consensus equilibria in \cite{IainD.Couzin2005Elad}. In some instances however, one may want a nonlinear model which always reaches consensus and the following result adapted from \cite{7330544} explains how to guarantee this.
\begin{proposition}
Suppose that for all edges, $f_{ij}$ has 0 as its unique root with $f_{ij}'(0)<0$. Then, the only equilibria of the system \cref{principal} in $\mathbb R ^n$ are the consensus states $\langle \mathbf 1 \rangle$. Moreover, all forward orbits converge to one of these states. 
\end{proposition}
\begin{proof}
For each edge, as $f_{ij}$ is odd, we can define $k_{ij}(|y|)=-f(y)/y$ for $y\not =0$ and $k_{ij}(0)=-f'(0)$. Then $f_{ij}(y)=-yk_{ij}(|y|)$ where  $k_{ij}$ is positive and continuous by the assumptions of this theorem. Now denote by $L(\mathbf x)$ the Laplacian matrix of the graph $\mathcal G$ with weights $k_{ij}(|x_i-x_j|)>0$. If we define,
\[E(\mathbf x)= \sum_{\{i,j\}\in \mathcal E}(x_i-x_j)^2,\]
one can show  that $\dot E(\mathbf x)= -\mathbf x^\top L(\mathbf x) \mathbf x$ (see \cite[Section IV]{7330544}) and as Laplacians are positive semi-definite $\dot E(\mathbf x)\leq 0$. As always it is enough to restrict our study to the state space $\langle \mathbf 1\rangle ^\bot$ where the only consensus state is the origin. Then, as $\mathcal G$ is connected we have $\dot E(\mathbf x)=0$ only at the origin, and  $E(\mathbf x)\rightarrow \infty$ if $||\mathbf x||\rightarrow \infty$. So applying   LaSalle's invariance principle \cite[Theorem 6.15]{book_Gerald_Teschl} to $E$, we conclude that all forward orbits in $\langle \mathbf 1\rangle ^\bot$ converge to the origin.
\end{proof}

\subsection{Stability of equilibrium points}
\label{sec:stability_eq}
To find the stability of an equilibrium point $\mathbf x^\star$ we consider the Jacobian matrix of \cref{principal} which we will denote by $J(\mathbf x^\star)$ or simply $J$. By the compact reformulation given in \cref{eq:compact_form2} and the chain rule we have,
\begin{equation}
    J(\mathbf x^\star)= d \hspace{3pt}D_{d^\top \mathbf x^\star} f\hspace{3pt} d ^\top
    \label{J = dDfd^T}.
\end{equation}
where $D_{\mathbf y} f$ denotes the Jacobian of $f$ at $\mathbf y$.  Note that $\mathbf 1$ is an eigenvector of eigenvalue 0 of $J$ and thus an equilibrium cannot be linearly stable in $\mathbb R^n$. However, it can be linearly stable in the dynamically invariant space  $\langle \mathbf 1\rangle ^\bot$, which is also invariant under the symmetric matrix $J$. From now on, linear stability will always be restricted to the space  $\langle \mathbf 1\rangle ^\bot$. Note that the eigenvalues of $J|_{\langle \mathbf 1\rangle ^\bot}$ are given by the ones in $\mathbb R^n$ while dropping one multiplicity of 0, which we will constantly use in our proofs. In particular, we have:



\begin{proposition}
Let $\mathbf x^\star$ be an equilibrium point of \cref{principal}. Then, if all eigenvalues of $J(\mathbf x^\star)$ (in $\langle \mathbf 1\rangle ^\bot$) are negative, $\mathbf x^\star $ is stable, whereas if $J(\mathbf x^\star)$ has a positive eigenvalue, $\mathbf x^\star $ is unstable. 
\label{prop:stability_equilibria}
\end{proposition}


\begin{remark}
In the original state space $\mathbb R^n$ we will never have attractive points. Indeed, any equilibrium point in $\langle \mathbf 1\rangle ^\bot$ spans a line of equilibrium points in the direction $\mathbf 1$ in  $\mathbb R^n$. However, if a point is stable in $\langle \mathbf 1 \rangle^\bot $ it will also be stable in $\mathbb R^n$ as all parallel planes have the same dynamics. 
\label{rem:stability}
\end{remark}

The conditions in the proposition above determine the linear stability of the equilibria in $\langle \mathbf 1 \rangle^\bot $. Thus, the stability will persist under small perturbations of the functions $f_{ij}$ as long as they remain odd. This will hold for all stability results, since they follow from \cref{prop:stability_equilibria}. 

As $J$ is a symmetric matrix, we can use the Courant minmax principle \cite{LotkerZvi2007Noda} which in particular implies,
 \begin{proposition}
 Let $M$ be a symmetric matrix and $\lambda _{\max} $ its maximum eigenvalue in an invariant subspace $V$. Then, for all $\mathbf x\in V$, 
 \[||\mathbf x||^2\lambda_{\max}\geq \mathbf x^\top M\mathbf x.\]
 Moreover, there exists a unitary vector $\mathbf z\in V$ such that, $\lambda _{\max} =\mathbf z^\top M\mathbf z$.
 \label{prop:minimax}
 \end{proposition}
 
In our context $V=\langle \mathbf 1 \rangle ^\bot$ and $M=J$. However, by the comment made above \cref{prop:stability_equilibria} if we find $\mathbf x\in \mathbb R^n$ such that $\mathbf x^\top J\mathbf x> 0$ it will be enough to deduce that $\lambda_{\max}>0$ in $\langle \mathbf 1 \rangle ^\bot$. Similarly, if we find $\mathbf x\in \mathbb R^n\setminus \langle \mathbf 1\rangle $ such that $\mathbf x^\top J\mathbf x \geq 0$ it will be enough to deduce that $\lambda_{\max}\geq 0$ in $\langle \mathbf 1 \rangle ^\bot$.

A direct consequence of this result are the following simple criterion for instability.  
\begin{proposition}
Let $\mathcal H$ be a cut-set of $\mathcal G$ and $\mathbf x^\star$ an equilibrium point of \cref{principal}. Then\footnote{Note that $f_{ij}'=f_{ji}'$ and is even, so the sum is well defined.}, 
\[\sum_{\{i,j\}\in \mathcal H}f'_{ij}(x^\star_i-x^\star_j)>0 \Longrightarrow \textrm{$\mathbf x^\star$ is unstable}. \]
In particular, for any $i\in \mathcal V$ we have, 
\[\sum_{j\sim i}f'_{ij}(x^\star_i-x^\star_j)>0 \Longrightarrow \textrm{$\mathbf x^\star$ is unstable}. \]
\label{prop:cut_set_unstability}
\end{proposition}
\begin{proof}
Let $\mathcal V=\mathcal V_1\cup \mathcal V_2$ be a partition of the nodes such that $\mathcal H$ is the associated cut-set and consider the vector $\mathbf x= \sum_{i\in \mathcal V_1} \mathbf e_i$. Now note that,
\begin{align*}
    \mathbf x ^\top J \mathbf x &= \sum_ {i\in \mathcal V_1} [J]_{i,i}+\sum_ {\substack{i,j\in \mathcal V_1\\j\sim i}} [J]_{i,j}= \sum_ {i\in \mathcal V_1} \sum _{j\sim i} f'_{ij}(x^\star_i-x^\star_j)-\sum_{i\in \mathcal V_1}\sum _{\substack{j\in \mathcal V_1\\j\sim i}}f'_{ij}(x^\star_i-x^\star_j)\\
    &=\sum _{i\in \mathcal V_1}\sum _{\substack{j\not \in \mathcal V_1\\j\sim i}}f'_{ij}(x^\star_i-x^\star_j) = \sum_{\{i,j\}\in \mathcal H}f'_{ij}(x^\star_i-x^\star_j)>0,
\end{align*}
and apply \cref{prop:minimax} together with \cref{prop:stability_equilibria}. 
\end{proof}

Informally the result above makes clear that if  $\mathbf x ^\star$ is stable, $f_{ij}'(x^\star_i-x^\star_j)\leq 0$ must hold for ``a substantial amount'' of edges, as we need to have at least one of them in each cut-set. Moreover, their module has to be big enough such that all corresponding sums are non-positive. We now proceed to show that ``a substantial amount'' can be expressed formally as containing a spanning tree. First we introduce some notation. 

Recall equation \cref{J = dDfd^T} and that $D_{d^\top (\mathbf x^\star)} f$ is a diagonal matrix, with values in the diagonal $f_{ij}'(x_i^\star -x_j^\star)$ for each $\{j,i\}\in {\mathcal E}$. Grouping the positive and negative values in two terms, we can decompose the Jacobian matrix as
\begin{equation}
J(\mathbf x^\star)=J = L^+- L^-,
     \label{J = L+ - L-}
\end{equation}
where $L^-$ (resp. $L^+$) is the weighted Laplacian of $\mathcal G^-$ (resp. $\mathcal G^+$) defined as the network $\mathcal G$ with weights given by $|f'_{ij}(x^\star_i-x^\star_j)|$ if $f'_{ij}(x^\star_i-x^\star_j)<0$ (resp. $f'_{ij}(x^\star_i-x^\star_j)>0$) and $0$ otherwise. Matrices of the form $L^+ -L^-$ are known as signed Laplacians and they do not retain many of the properties of standard Laplacians, for instance they may have both positive and negative eigenvalues. A general study on the signs of their spectrum
can be found in \cite{BronskiJaredC.2014STfD}, where an alternative proof of the following result is given. 

\begin{proposition}
Given an equilibrium point $\mathbf x^\star$ of \cref{principal} we have, 
\begin{itemize}
    \item if $L^+=\mathbf 0_{n,n}$ and $\mathcal G^-$ is connected,  $\mathbf x^\star $ is stable;
    \item if $\mathcal G^-$ is disconnected then $\mathbf x^\star $ is not linearly stable (in $\langle \mathbf 1\rangle ^\bot$). Moreover, if there is an edge in $\mathcal G^+$ between two distinct connected components of $\mathcal G^-$, $\mathbf x^\star$ is unstable.
\end{itemize}
\label{prop:G-connected}
\end{proposition}
\begin{proof}
See \cref{app:connected}.

\end{proof}


This proposition is useful to study some simple situations. For instance, if a network has an edge $\{i,j\}$ with no cycles going trough it, e.g. the connecting edge of the barbell graph, and $f'_{ij}(x^\star_i-x^\star_j)>0$ then $\mathbf x^\star$ is unstable. To be able to deal with more complex situations we will introduce the concept of effective resistances in \cref{sec:effective_resistance}. 

%

\subsubsection{Schur complement reduction}

In this section we show how to use the Schur complement to reduce the dimensionality of the stability problem. Essentially we will be able to remove all nodes which are not incident to any edge in $\mathcal G^+$. Clearly, this technique is very powerful when  $\mathcal G^+$ contains few edges, see \cref{thm:resistencia:one_edge} for the case with a single edge, and \cref{app:Schur} for a concrete example. 

We start by introducing the Schur complement for a symmetric matrix.

\begin{definition}
 Given a $n\times n$ symmetric matrix $M$, and $\mathcal N, \mathcal  M \subset \{1,\dots,n\}$, we denote by  $M_{\mathcal N\mathcal M}$ the submatrix of $M$ with rows indexed by $\mathcal N$ and columns by $\mathcal M$. Then, if $M_{\mathcal N\mathcal N}$ is invertible,  we define the \emph{Schur complement} of $M$ respect to $\mathcal N$ by,
 \begin{equation}M/\mathcal N:=M_{\mathcal N^c \mathcal N^c}-M_{\mathcal N^c \mathcal N} M_{\mathcal N \mathcal N }^{-1}\mathcal M_{\mathcal N \mathcal N^c}.
 \label{eq:schur_def}\end{equation}
\end{definition}

The Schur complement has many interesting properties and applications as reviewed in  \cite{Schur_complement}. 
In our context, its two key features are that $\In (M)= \In (M/\mathcal N)+\In ( M_{\mathcal N  \mathcal N})$,
see \cite[Theorem 1.6]{Schur_complement}, and that if $M$ is a weighted Laplacian then $M/\mathcal N$ is also a Laplacian \cite{DorflerFlorian2018ENaA}. 

Recall that finding the stability of an equilibrium $\mathbf x^\star$ reduces to the study of the inertia of $J$. Here it will be convenient to assume $\mathcal G^-$ connected and one can then use \cref{prop:G-connected} to deduce the stability for the general case. Denote by $\mathcal N$ the set of nodes  non-incident to any edge in $\mathcal G^+$ and assume that  $\mathcal N$ and $\mathcal N^c$ are non-empty. Then, reordering the nodes we have, 
\begin{equation}J=\begin{pmatrix}L^+_{\mathcal N^c\mathcal N^c}-L^-_{\mathcal N^c\mathcal N^c} & -L^-_{\mathcal N^c\mathcal N}\\ -L^-_{\mathcal N\mathcal N^c}&-L^-_{\mathcal N\mathcal N}\end{pmatrix}.\label{eq:J_schur}\end{equation}
As $L^-_{\mathcal N\mathcal N}$ is a principal submatrix of a connected Laplacian, it is positive definite, so invertible. Thus, we can consider the Schur complement $J/\mathcal N$ and get,
\[\In (J)= \In \left (J/\mathcal N\right )+\In \left  ( J_{\mathcal N  \mathcal N} \right )=\In \left (J/\mathcal N\right )+\In \left  (- L^-_{\mathcal N  \mathcal N} \right ).\]
As $- L^-_{\mathcal N  \mathcal N}$ is negative definite, the stability of $\mathbf x^\star$ is determined by $\In \left (J/\mathcal N\right )$. Moreover, from \cref{eq:J_schur} it is clear that  
\begin{equation}
J/\mathcal N =L^+_{\mathcal N^c \mathcal N ^c}- L^-/\mathcal N.
\label{eq:J/N}
\end{equation}
where $L^-/\mathcal N$ is a Laplacian matrix, as it is the Schur complement of a Laplacian. Note that $L^+_{\mathcal N^c \mathcal N ^c}$ is also a Laplacian matrix, as by definition of $\mathcal N$ all other entries of $L^+$ are zeros. Thus, expression \cref{eq:J/N} is analogous to \cref{J = L+ - L-} but in dimension $|\mathcal V|-|\mathcal N|$, and we can interpret it as the stability problem for a smaller network. 
In summary we have shown:

\begin{theorem}
Let $\mathbf x^\star$ be an equilibrium point of \cref{principal}, assume that $\mathcal G^-$ is connected and let $\mathcal N$ be the nodes not incident to any edge in $\mathcal G^+$. Then, the  stability of $\mathbf x^\star$  is given by the linear stability of $J/\mathcal N$ in $\langle \mathbf 1\rangle ^\bot$. 
\label{thm:schur_comp}
\end{theorem}

This theorem will only be useful when $|\mathcal N|$ is large, in particular if $\mathcal G^+$ is connected, then $|\mathcal N|=0$, and we do not get any reduction. One small detail we have ignored so far is that in  \cref{eq:J/N} one may have edges that are contained in both Laplacians $L^+_{\mathcal N^c \mathcal N ^c}$ and $ L^-/\mathcal N$, in contrast to \cref{J = L+ - L-}. This can be used to apply the Schur complement to a new decomposition and reduce the dimensions even further, see \cref{app:Schur}.


\subsubsection{Effective resistance criteria}
\label{sec:effective_resistance}
In this section we will apply ideas from electrical networks, and in particular the concept of effective resistance, to determine the stability of the equilibrium points of \cref{principal}. This application of effective resistance was introduced in the recent work by Devriendt et al. \cite{karel2020arXiv} for a system of the form \cref{principal}, and by Tyloo et al. \cite{Tyloo_resistances_in_Kuramoto} for Kuramoto-like systems. While the analysis in \cite{karel2020arXiv} was based on unweighted networks, we will consider weights on $\mathcal{G}$ in our more general setting, leading to improved results. 
We will show that when a pair of nodes are better connected\footnote{Better connected in the sense of effective resistance, i.e. the effective resistance between them in $\mathcal G^+$ is smaller than the one in $\mathcal G^-$.}  in $\mathcal G^+$ than in $\mathcal G^-$ then $\mathbf x^\star $ is unstable. Moreover, when $\mathcal G^+$ contains a single edge, this condition is tight. 

\begin{definition}
 The \emph{effective resistance} $r_{ij}$ between a pair of nodes $i$ and $j$ in a weighted network with Laplacian $L$, is defined as (see, e.g. \cite{DorflerFlorian2018ENaA,KleinD.1993Rd})
\begin{equation*}
    r_{ij}:= (\mathbf e_i-\mathbf e_j)^\top L^\dag (\mathbf e_i-\mathbf e_j)
\end{equation*}
if $i$ and $j$ are in the same connected component, and $r_{ij}=\infty$ otherwise. 
\label{def:ef.res.}
\end{definition}

Note that the effective resistance between two adjacent nodes, with no other path between them, is given by the inverse of the edge weight.  Moreover, this definition satisfies the well known properties of resistors in series and in parallel, see \cref{supp:sec:parallel}. A generalization of this fact is given by the invariance of effective resistance by the Schur complement, also known as Kron reduction \cite{DorflerFlorian2018ENaA}. That is, the effective resistance between two nodes $i,j\not \in \mathcal N\subset \mathcal V$ is the same in $L$ and in $L/\mathcal N$.

 Recall that $J=L^+-L^-$  and denote by $r_{ij}^+$, $r_{ij}^-$ the effective resistance between $i$ and $j$ of the corresponding networks $\mathcal G^+$, $\mathcal G^-$.
With this notation we are prepared to state our first result, which gives tight conditions for the stability when $\mathcal G^+$ has a single edge.
 \begin{theorem}
Let $\mathbf x^\star$ be an equilibrium point of \cref{principal} such that there exists a unique edge $\{i,j\}$  with $f'_{ij}(x^\star_i-x^\star_j)\geq 0$ and assume that $f'_{ij}(x^\star_i-x^\star_j)>0$. Then, 
\begin{align*} 
f'_{ij}(x^\star_i-x^\star_j)r_{ij}^-<1 &\Longrightarrow \textrm{$\mathbf x^\star$ is stable}; \\ 
f'_{ij}(x^\star_i-x^\star_j)r_{ij}^->1 &\Longrightarrow \textrm{$\mathbf x^\star$ is unstable} .
\end{align*}
\label{thm:resistencia:one_edge}
\end{theorem}
\begin{proof}
If $\mathcal G^-$ is disconnected we have $r_{ij}^-=\infty$  and $\mathbf x^\star$ is unstable due to  \cref{prop:G-connected}.  If $\mathcal G^-$ is connected, by  \cref{thm:schur_comp} the stability is determined by $J/\mathcal N$. As $\mathcal N^c=\{i,j\}$, we have, 
\[J/\mathcal N = L^+_{\{i,j\},\{i,j\}} - L ^-/\mathcal N= \begin{pmatrix}f'_{ij}(x^\star_i-x^\star_j) & -f'_{ij}(x^\star_i-x^\star_j)\\-f'_{ij}(x^\star_i-x^\star_j) & f'_{ij}(x^\star_i-x^\star_j)\end{pmatrix}-\begin{pmatrix}w_{ij}^- & -w_{ij}^-\\-w_{ij}^- & w_{ij}^-\end{pmatrix},\]
for certain weight $w_{ij}^-$. Thus, $\mathbf x^\star$ is stable if $f'_{ij}(x^\star_i-x^\star_j)-w_{ij}^-<0$ and unstable if $f'_{ij}(x^\star_i-x^\star_j)-w_{ij}^->0$. Finally, we use that effective resistance is invariant under Schur complement, so  $r_{ij}^-=1/w_{ij}^-$.
\end{proof}

\begin{remark}
If there is a unique edge with $f'_{ij}(x^\star_i-x^\star_j)>0$ and $\mathcal G^-\cup \{i,j\}$ is connected, the result holds with the same arguments. 
\label{rem:one_edge}
\end{remark}

We can take an analogous approach when $\mathcal G^+$ contains two edges with a common node, or in fact any particular small configuration. Then, using Sylvester's criterion \cite{Semidef}, one gets a couple of inequalities on the entries of $J/\mathcal N$ (which can be related to effective resistance) that determine the stability. This operation gets quite messy so, instead,
 we expand by brute force \cref{thm:resistencia:one_edge} to the general case, but we lose the tightness of the bounds.


\begin{proposition}
Let $\mathbf x^\star$ be an equilibrium point of \cref{principal}. Then, if $\mathcal G^+\cup \mathcal G^-$ is connected,
\begin{align*} 
\sum_{\{i,j\}\in \mathcal G^+ }f'_{ij}(x^\star_i-x^\star_j)r_{ij}^-<1 &\Longrightarrow \textrm{$\mathbf x^\star$ is stable}; \\ 
\exists \{i,j\}\in \mathcal G^+ : f'_{ij}(x^\star_i-x^\star_j)r_{ij}^->1 &\Longrightarrow \textrm{$\mathbf x^\star$ is unstable} .
\end{align*}
\label{thm:resistencia:mult_edge}
\end{proposition}
\begin{proof}
See \cref{app:resistance_thm}.
\end{proof}

Both \cref{thm:resistencia:one_edge} and \cref{thm:resistencia:mult_edge} are adapted results from \cite{karel2020arXiv} where a particular coupling function is considered. We have opted to give completely different proofs, but  it is possible to generalize the arguments from \cite{karel2020arXiv} to our context with some work.
Let us now present another approach which leads to a tighter instability condition. First, we need the following result, see \cref{app:lem:unstable} for the proof.



\begin{lemma}
Let $A, B$ be symmetric positive semi-definite matrices of the same dimension. Then\footnote{All maximums are taken over non-zero vectors.}, 
\[\max_{\mathbf x\bot \ker B}\frac{\mathbf x^\top  A\mathbf x}{\mathbf x^\top B\mathbf x}= \max_{\mathbf x\bot \ker A}\frac{\mathbf x^\top  B^\dag \mathbf x}{\mathbf x^\top A^\dag \mathbf x}.\]
\label{lem:instabiliy}
\end{lemma}


In \cref{supp:sec:stability B^T/A^T} we discuss how this result could also be used to get stability conditions.


\begin{theorem}
Let $\mathbf x^\star$ be an equilibrium point of our system and $i,j$ a pair of nodes. Then,
\[r^-_{ij}>r^+_{ij} \Longrightarrow \textrm{$\mathbf x^\star$ is unstable}.\] 
\label{thm:unstability_r^-_ij r^+_ij}
\end{theorem}
\begin{proof}
Denote by $\lambda_{\max}$ the maximum eigenvalue of $J$ in $\langle \mathbf 1 \rangle^\bot$, then by \cref{prop:minimax} we have, 
\[\lambda_{\max}\geq \max_{\substack{\mathbf x\bot \ker L^-\\||\mathbf x||=1}} {\mathbf x}^\top L^+\mathbf x- {\mathbf x}^\top L^-\mathbf x=\max_{\substack{\mathbf x\bot \ker L^-\\||\mathbf x||=1}}{\mathbf x}^\top L^-\mathbf x\left (\frac{{\mathbf x}^\top L^+\mathbf x}{\mathbf x^\top L^-\mathbf x} -1 \right)\]
and as $L^-$ is positive semi-definite,
\begin{equation}
    \textrm{sign}(\lambda_{\max})\geq \textrm{sign}\left (\max_{\mathbf x\bot \ker L^-} \frac{\mathbf x^\top L^+\mathbf x}{\mathbf x^\top L^-\mathbf x} -1 \right).
    \label{eq:sign_lambda}
\end{equation}
Now, if $r_{ij}^+ =\infty $ the condition of this theorem is never satisfied and there is nothing to prove. If $r_{ij}^-=\infty $ and $r_{ij}^+<\infty $, $\mathbf x^\star $ is unstable by \cref{prop:G-connected}. So we may assume that $i$ and $j$ are in the same connected component of $\mathcal G^+$ and $\mathcal G^-$. Then, $(\mathbf e_i-\mathbf e_j)\bot  \ker L^+$ and by the previous lemma, 
\[\max_{\mathbf x\bot \ker L^-}\frac{{\mathbf x}^\top  L^+ {\mathbf x}}{{\mathbf x}^\top L^- {\mathbf x}}= \max_{{\mathbf x}\bot \ker L^+}\frac{{\mathbf x}^\top  (L^-)^\dag {\mathbf x}}{{\mathbf x}^\top (L^+)^\dag {\mathbf x}}\geq \frac{(\mathbf e_i-\mathbf e_j)^\top  (L^-)^\dag (\mathbf e_i-\mathbf e_j)}{(\mathbf e_i-\mathbf e_j)^\top (L ^+)^\dag (\mathbf e_i-\mathbf e_j)}=\frac{r_{ij}^-}{r_{ij}^+}.\]
So if $r^-_{ij}>r^+_{ij}$, then  $\frac{r_{ij}^-}{r_{ij}^+}-1>0$, and thus by \cref{eq:sign_lambda}, $\lambda_{\max}>0$. 
\end{proof}

\begin{remark}
When there is a unique edge  $\{i,j\}$ with positive derivative then, $r^+_{ij}=1/f'(x^\star_i-x^\star_j)$ and the instability condition coincides with \cref{thm:resistencia:one_edge}. In general $r_{ij}^+\leq 1/f'(x^\star_i-x^\star_j)$, and the new instability condition is tighter than \cref{thm:resistencia:mult_edge}.  
\end{remark}

It has been shown that the effective resistance is a distance function on the nodes of a network which in some sense encapsulates how well connected the different pairs of nodes are \cite{KleinD.1993Rd,DorflerFlorian2018ENaA}. With this point of view, the theorem above asserts that if there exist a pair of nodes which are better connected in $\mathcal G^+$ than in $\mathcal G^-$, then we have instability. For the case of a unique non-negative edge, \cref{thm:resistencia:one_edge} asserts that this condition is tight, and only needs to be checked for the non-negative edge. \Cref{thm:unstability_r^-_ij r^+_ij} does not give tight conditions in general, but they might be tight for certain types of networks (see \cref{supp:sec:instability_kn}).


\subsection{Union of two networks in a node}
\label{section:union_graphs}
Assume that $\mathcal G$ contains a \emph{cut-node}, i.e. a node that when ``deleted'' increases the number of connected components of $\mathcal G$. Equivalently $\mathcal G$ is the union of two subnetworks ${\mathcal A}$, ${\mathcal B}$ that intersect in a single node, also known as a \emph{coalescence} of $\mathcal{A}$ and $\mathcal{B}$ \cite{Coalescence_dynamics}. In this section we will show that the equilibrium points of $\mathcal G$ are exactly the combination of equilibrium points of ${\mathcal A}$ and ${\mathcal B}$. Moreover, we will show that the linear stability of an equilibrium point in $\mathcal G$ is determined by the stability of the corresponding equilibrium points in ${\mathcal A}$ and ${\mathcal B}$. 
These results are quite surprising, as there are several ways to coalesce two networks, which can produce significantly different structures. Our results show that from a basic dynamical point of view, the resulting networks will be equivalent. 
A comparison between the dynamical properties of networks and their coalescence was also considered in \cite{Coalescence_dynamics} for a system of coupled oscillators, where coalescence was shown to lead to an increased time to reach synchrony.

Without loss of generality we may assume that the nodes of  ${\mathcal A}$ and  ${\mathcal B}$ are $\{1,\dots,\tilde n\}$ and $\{\tilde n,\dots, n\}$ respectively, and let $n'=n-\tilde n+1$. As all edges in $\mathcal G$ are in ${\mathcal A}$ or ${\mathcal B} $ exclusively, the edge state of $\mathcal G$ is given by the edge states of  ${\mathcal A}$ and ${\mathcal B}$ combined. That is, reordering the components if necessary $\mathbf y =(\mathbf y_{\mathcal A},\mathbf y_{\mathcal B})$. Note that with these assumptions $f(\mathbf y)=(f_{\mathcal A}(\mathbf y_{\mathcal A}),f_{\mathcal B}(\mathbf y_{\mathcal B}))$ and the generating cycles of $\ker d$ are entirely contained in ${\mathcal A}$ or ${\mathcal B}$. Then, by \cref{eq:equlibrium_edges} it is straightforward to show that,
\[\mathscr E_{\mathbf y}= \iota_{\mathcal A}(\mathscr E_{\mathbf y} ^{\mathcal A})\oplus \iota_{\mathcal B}(\mathscr E_{\mathbf y} ^{\mathcal B}),\]
where, $\mathscr E_{\mathbf y}^{\mathcal A}$ are the equilibria on $\mathcal A$ and   $\iota_{\mathcal A} (\mathbf y_{\mathcal A})=(\mathbf y_{\mathcal A}, \mathbf 0_{\mathcal B}) $ (resp. for $\mathcal B$), see \cref{supp:sec:iota_A} for a detailed exposition. Hence,  $\mathbf y=(\mathbf y_{\mathcal A},\mathbf y_{\mathcal B})\in \mathscr E_{\mathbf y}$ if and only if $\mathbf y_{\mathcal A}\in \mathscr E_{\mathbf y}^{\mathcal A}$ and $\mathbf y_{\mathcal B}\in \mathscr E_{\mathbf y}^{\mathcal B}$. 

Now we move on to prove that the stability can  also be studied by restricting ourselves to ${\mathcal A}$ and ${\mathcal B}$. We need the following technical result,  proved in \cref{app:union_graphs}. 

\begin{lemma}
 Let $\mathcal G$ be a network with $m$ edges, $d$ its boundary map and $\tilde d$ the same matrix with a row deleted. Let $M\in \mathbb R^{m\times m}$ be a diagonal matrix. Then, the number of positive (resp. negative) eigenvalues of $dMd^\top $ and $\tilde dM\tilde d^\top$ coincide. 
 \label{union_graph:lema}
\end{lemma}

 Recall from \cref{prop:stability_equilibria} that the stability only depends on the sign of the eigenvalues of $J$. We have, 
\[J=d D_{\mathbf y} f d^\top=\begin{pmatrix}
 d_{\mathcal A}  D_{\mathbf y_{\mathcal A}} f_{\mathcal A} d_{\mathcal A}^\top & \mathbf 0_{\tilde n,n'-1} \\
\mathbf 0_{n'-1,\tilde n} &\mathbf 0_{n'-1,n'-1}
\end{pmatrix}+\begin{pmatrix}
 \mathbf 0_{\tilde n-1,\tilde n-1}& \mathbf 0_{\tilde n-1,n'}\\
\mathbf 0_{n',\tilde n-1} & d_{\mathcal B} D_{\mathbf y_{\mathcal B}} f_{\mathcal B} d_{\mathcal B}^\top
\end{pmatrix},\]
where in the $\tilde n$th row/column we can have non-zero entries in both matrices. Now if we denote by $\tilde d$, $\tilde d_{\mathcal A}$, $\tilde d_{\mathcal B}$ the respective matrices when deleting the row corresponding to the node $\tilde n$, we find
\[\tilde d D_{\mathbf y} f \tilde d^\top=\begin{pmatrix}
 \tilde d_{\mathcal A}  D_{\mathbf y_{\mathcal A}} f_{\mathcal A} \tilde d_{\mathcal A}^\top & \mathbf 0_{\tilde n-1,n'-1}\\
\mathbf 0_{n'-1,\tilde n-1} & \tilde d_{\mathcal B} D_{\mathbf y_{\mathcal B}} f_{\mathcal B} \tilde d_{\mathcal B}^\top
\end{pmatrix}.\]
So the spectrum of $\tilde d D_{\mathbf y} f \tilde d^\top$ is the union of the spectrum of $\tilde d_{\mathcal A}  D_{\mathbf y_{\mathcal A}} f_{\mathcal A} \tilde d_{\mathcal A}^\top$ and 
$ \tilde d_{\mathcal B} D_{\mathbf y_{\mathcal B}} f_{\mathcal B} \tilde d_{\mathcal B}^\top$. Thus, by \cref{union_graph:lema}, the signs of the  spectrum of $ d D_{\mathbf y} f d^\top$  are the union of the signs in the spectrum of  $ d_{\mathcal A}  D_{\mathbf y_{\mathcal A}} f_{\mathcal A} d_{\mathcal A}^\top$ and $ d_{\mathcal B}  D_{\mathbf y_ {\mathcal B}} f_{\mathcal B} d_{\mathcal B}^\top$ (with a zero removed due to dimensionality). In particular,  the linear stability of $\mathbf y$ is determined by the linear stability of $\mathbf y_{\mathcal A}$ and $\mathbf y_{\mathcal B}$.


\subsubsection{Tree of motifs}
\label{sec:tree_motifs}
\begin{figure}[htbp]
    \centering
        {\includegraphics[clip, trim=0cm 0cm 0cm 0cm, width=0.6\textwidth]{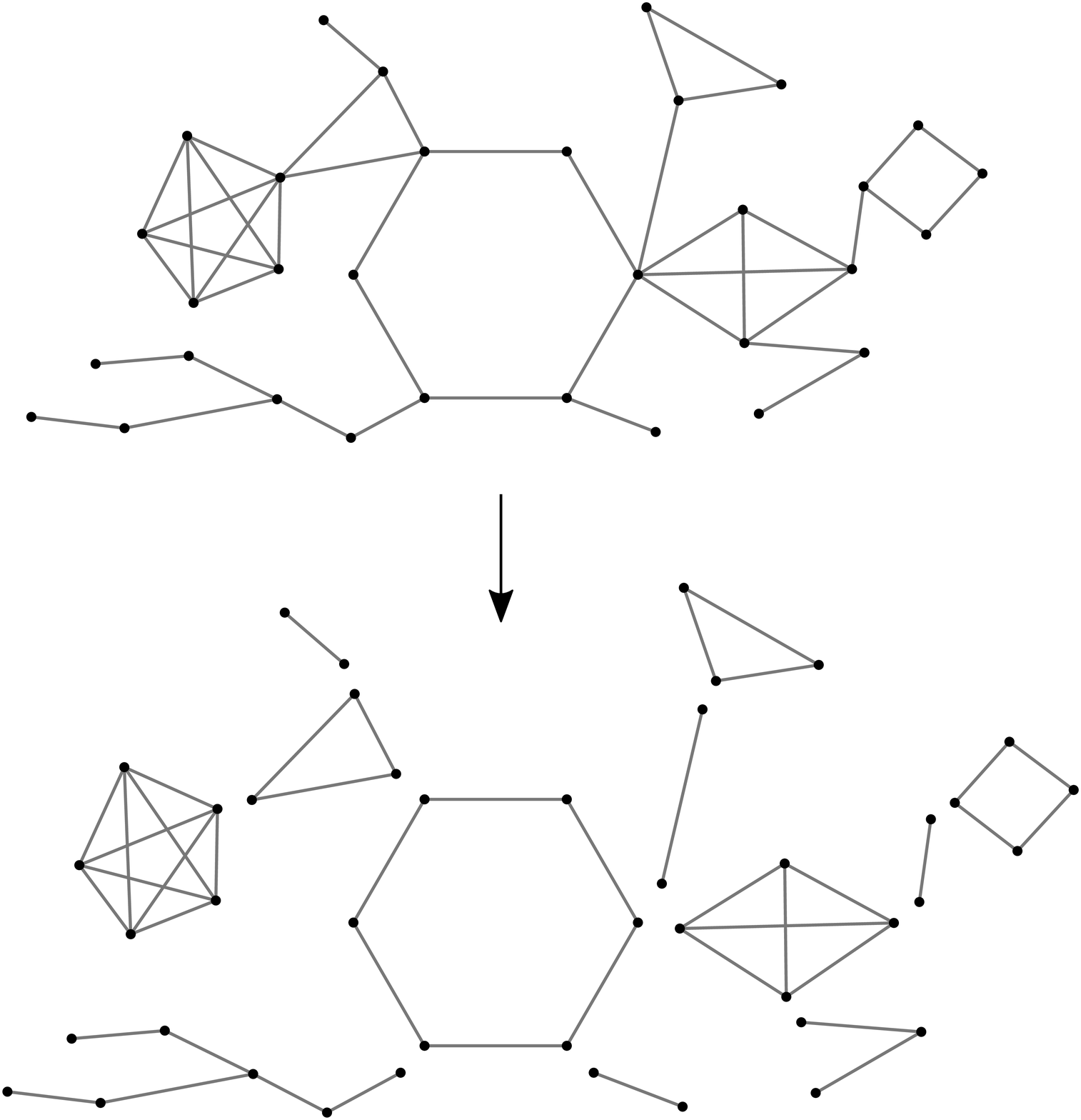}}
    \caption{In the upper half a tree of motifs is shown where the motifs are trees, cycle and complete graphs. In the lower half, the motifs are shown separated from each other. These are the subnetworks that we would need to understand dynamically to comprehend the behaviour of the whole network. }
    \label{fig:geo_graph}
\end{figure}

Starting from a collection of ``small'' networks which are understood dynamically\footnote{That is, we know their equilibrium points and their linear stability.}, we can create new networks by recursively joining (coalescing) these smaller networks in single nodes. By our results above, the dynamics of this new network will be equally well understood by simply composing the properties of the subsystems. We refer to these smaller networks as \emph{motifs} as they play the role of dynamical and structural building blocks, and refer to the larger network as a \emph{tree of motifs} as they consist of motifs interconnected in a tree (loopless) way, see \cref{fig:geo_graph}.

Equivalently our result can be stated as: given an arbitrary network the task of finding equilibria and their stability can be reduced to the same task for each of its \emph{blocks}\footnote{A block of a network is a maximal connected subnetwork that does not contain any cut-nodes  of itself.}.
While some networks will consist of a single block, for many other networks this procedure will drastically reduce the dimensions in which we are working. In \cite{Kiss2015} similar results are shown for Markovian $SIR$ epidemic dynamics. 


 \section{Applications to concrete networks}
 \label{chapter:application_to_network}
In this section we will focus on studying the equilibrium points of \cref{principal} and their stability in some specific types of networks. Concretely we will study tree, cycle and complete graphs. 
Although it may seem that these families of networks are quite specific, using the main result from \cref{section:union_graphs},  we will be able to tackle trees of motifs for the motifs: trees, cycles and complete graphs. See \cref{fig:geo_graph} for an example and  \cref{supp:sec:example} for explicit computations.
 Such trees of motifs appear as the construct of certain growing graph models (see \cite[Section  3.5.2]{aynaud2013multilevel} and \cite[Section 2]{growing_graph_model}), and the special case where all motifs are complete graphs is also known as block graphs \cite{block_graph}.


\subsection{Tree graphs}
In this section we will study the case when $\mathcal G$ has the  simplest topology, which in our context means that it does not contain any cycles, i.e.  $\mathcal G$ is a tree. We will show explicit expressions of the equilibrium points of \cref{principal} given the roots of the functions $f_{ij}$, and simple criterion for their stability. These results will directly follow from the previous section, as tree graphs are in particular trees of motifs, where all motifs are simply a pair of connected nodes. 

\begin{proposition}
In a tree graph we have, 
\[\mathscr E_{\mathbf y}=\tilde{\mathscr E}_{\mathbf y}=f^{-1}(\mathbf 0)=\left \{ (y_{\{i,j\}})_{\{i,j\}\in \mathcal E}\in \mathbb R^m \textrm{ \textnormal{:} } f_{ij}(y_{\{i,j\}})=0 \textrm{ \textnormal{for all} } \{i,j\}\in \mathcal E\right \},\]
Moreover, the signs of the eigenvalues of $J(\mathbf x^\star)$ in $\langle \mathbf 1 \rangle^\bot$ are given by $\textrm{\textnormal{sign}}(f'_{ij}(x^\star_i-x^\star_j))$ for $\{i,j\}\in \mathcal E$. In particular, if there exists $\{i,j\}\in \mathcal E$ such that $f'_{ij}(x^\star_i-x^\star_j)>0$, then $\mathbf x^\star$ is unstable, whereas $\mathbf x^\star$ is stable if $f'_{ij}(x^\star_i-x^\star_j)<0$ for all $\{i,j\}\in \mathcal E$.
\label{prop:tree_stability}
\end{proposition}
\begin{proof}
First, as trees have no cycles, $\ker d=\mathbf 0$ so from \cref{eq:equlibrium_edges} and the definition of $\tilde{\mathscr E}_{\mathbf y}$ we get  $\mathscr E_{\mathbf y}=\tilde{\mathscr E}_{\mathbf y}=f^{-1}(\mathbf 0)$. For the signs of the eigenvalues apply \cref{section:union_graphs} inductively on edges.
\end{proof}

In particular, in a tree, the set of equilibria (on the edge space) and their linear stability only depends on the type of coupling functions $f_{ij}$, and not on the specific arrangement of the edges.  For instance, the equilibria and stability in a line graph  will be the same as in a star graph of the same size.  Note that the result above not only gives us the stability of the equilibria, but also finds the signs of all eigenvalues of $J$. This will be useful for the following result, which informally asserts that any random initial condition will converge to one of the stable equilibrium points or  ``infinity''. 

\begin{proposition}
Let $\mathcal G$ be a tree  where for all edges the roots of $f_{ij}$ are simple. 
Then, for Lebesgue almost every point $\mathbf x\in \langle \mathbf 1\rangle ^\bot$  its forward orbit converges to a stable equilibrium or is unbounded.
Moreover, the stable equilibrium points are given by,
\[ \left \{\mathbf x\in \langle \mathbf 1\rangle^\bot \textrm{\textnormal{ : }}  f_{ij}(x_i-x_j)=0 \textrm{\textnormal{ and }} f_{ij}'(x_i-x_j)< 0 \textrm{\textnormal{ for all }} \{i,j\}\in \mathcal E\right \}.\]
\label{prop:mesure0}
\end{proposition}
\begin{proof}
By the previous proposition, all equilibrium points are hyperbolic, thus isolated,
 and the stable ones are given by the expression above. Then,  in gradient dynamics, LaSalle's invariance principle \cite[Theorem 6.15]{book_Gerald_Teschl} guarantees that all forward bounded orbits converge to a set of equilibrium points. As they are isolated, orbits converge to a single point. It also follows that there are countably many equilibrium points so it is enough to prove that the set of points converging to an unstable equilibrium has null measure. This holds as its stable manifold has at least codimension 1 (see  \cref{lemma:grad_mesure0} for details). 
\end{proof}

\begin{remark}
In the previous result, if for all edges, $f_{ij}$ has 3 simple roots with $f'_{ij}(0)<0$, then the only stable equilibrium is the origin. So almost every forward bounded orbit converges to it. 
\end{remark}


Imposing more conditions on $f_{ij}$,  one can show that all forward orbits are bounded \cref{prop:descomp_state_W}. To get a flavour of how a complete description of the dynamics of \eqref{principal} may look, we give the phase portrait of trees in low dimension, see \cref{sec:pahse_portrait}. 

\subsection{Cycle graphs}
\label{sec:cycle_graph}
Let $\mathcal C_n$ be the cycle graph  with $n$ nodes, i.e a line graph where the end nodes have been joined, and let  $\mathbf{x}^\star$ be an equilibrium point of \cref{principal} on this network.
Then if $f_{ij}'(x_i^\star-x^\star_j)<0$ for all edges, $\mathbf{x}^\star$ is stable by \cref{prop:G-connected}.  If instead there exist two edges with non-negative derivatives and one of them is positive, $\mathbf x^\star$ is unstable, by \cref{prop:G-connected}. Assume now that there is a unique edge $\{i,j\}$ with non-negative derivative, and moreover $f_{ij}'(x_i^\star-x^\star_j)\not =0$. Then, by  \cref{thm:resistencia:one_edge} we have,
\begin{align*} 
f'_{ij}(x^\star_i-x^\star_j) \cdot\sum_{\{i',j'\}\in \mathcal E\setminus\{\{i,j\}\}}\frac{1}{|f'_{i'j'}(x_{i'}^\star-x_{j'}^\star)|}<1 &\Longrightarrow \textrm{$\mathbf x^\star$ is stable}; \\ 
f'_{ij}(x^\star_i-x^\star_j)\cdot\sum_{\{i',j'\}\in \mathcal E\setminus\{\{i,j\}\}}\frac{1}{|f'_{i'j'}(x_{i'}^\star-x_{j'}^\star)|}>1 &\Longrightarrow \textrm{$\mathbf x^\star$ is unstable},
\end{align*}
where we have used the formula for effective resistances in series, see  \cref{supp:eq:series_res}. Thus, we have tight conditions for the stability of $\mathbf x^\star$ except for the case when $\max_{\{i,j\}\in \mathcal E} f_{ij}'(x_i^\star-x_j^\star)=0$ or when there is a unique edge with non-negative derivative and the expression above equals 1.

\begin{figure}[htbp]
    \centering
        {\includegraphics[clip, trim=0cm 0.5cm 0cm 0.5cm, width=0.4\textwidth]{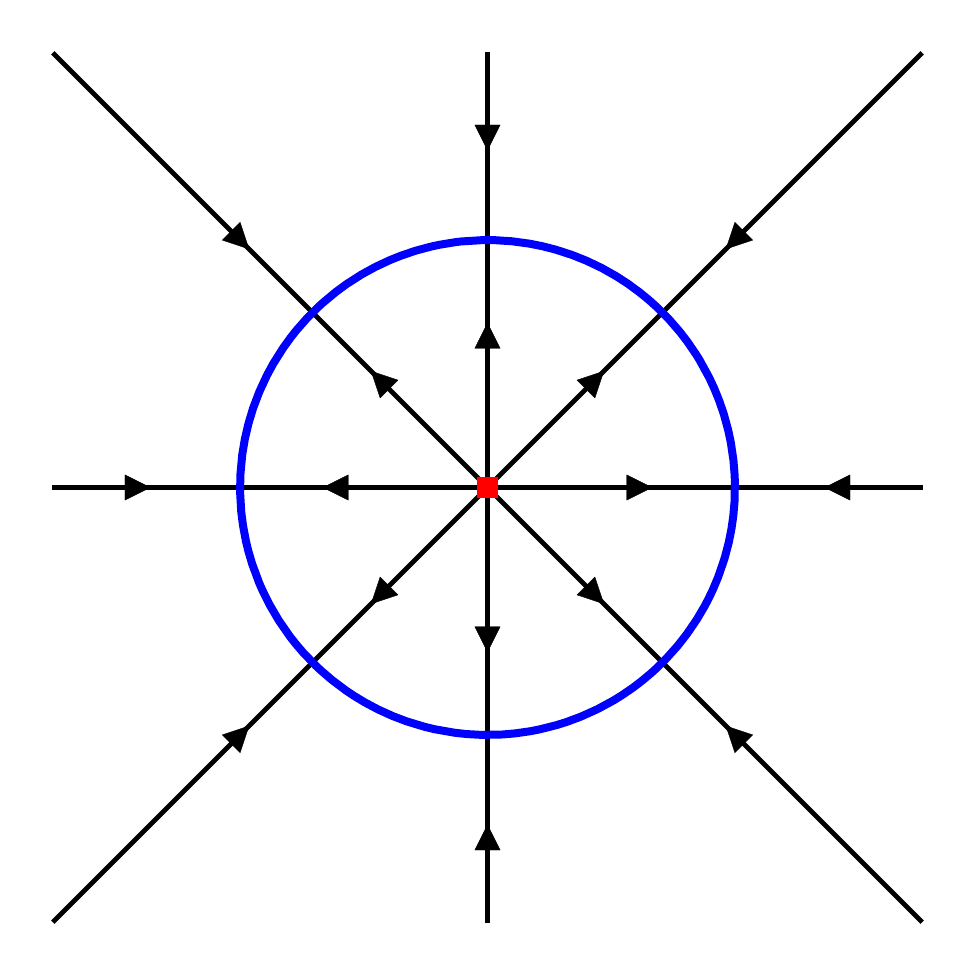}}
    \caption{Phase portrait of \cref{principal} for $\mathcal C_3$ and $\hat f(y)=\lambda y -y^3$ with $\lambda>0$ in a plane perpendicular to $\mathbf 1$. A circle of stable equilibrium points of radius $\sqrt{2\lambda/3}$ is depicted in \emph{blue}, and a source equilibrium is shown by a \emph{red square}. All other orbits are contained in lines going through the source point.}
    \label{fig:circle}
\end{figure}

In contrast with trees, where $\tilde {\mathscr E}_{\mathbf x}=\mathscr E_{\mathbf x}$, this equality will not always hold for cycle graphs (which may be considered the next simplest topology). This was already identified in \cite{Onbif} for $\mathcal C_3$ taking  $f_{ij}(y)=\hat f(y)=\lambda y-y^3$ with $\lambda>0$ for all $\{i,j\}\in \mathcal E$. They noted that in the original node space $\mathbb R^3$, the subspace $\langle \mathbf 1\rangle $ consists of       unstable equilibrium points, and that there is a cylinder with axis $\langle \mathbf 1\rangle $, consisting of stable equilibria.  We note that, 
\[H(x_1,x_2,x_3)=\frac{x_1-x_2}{x_1-x_3},\]
is a first integral  of this system, so that the phase portrait in the parallel planes to $\langle\mathbf 1 \rangle^\bot $ is as depicted in \cref{fig:circle}.  Notice that we get a continuum of equilibrium points. 

We now study when and why this phenomenon of a curve of equilibria in the node space $\langle \mathbf 1\rangle^\bot$ occurs  for a general function $\hat f$. As  $n=m$ and following the comments made in  \cref{sec:equilibrium_points}, we suspect that this set is at most one dimensional. In this context it will be convenient to work with a slight modification of the edge coordinates,
\[z_1=y_{\{1,2\}},\hspace{5pt}\dots\hspace{3pt},\hspace{5pt} z_{n-1}=y_{\{n-1,n\}},\hspace{5pt}z_{n}=-y_{\{1,n\}}.\]
Recall that $\ker d$ is generated by the cycles so in $\mathbf z$ coordinates, $\ker d= \langle \mathbf 1\rangle$. Thus from \cref{eq:equlibrium_edges} we get,
\begin{equation}
    \mathscr E_{\mathbf z}= \langle {\mathbf 1} \rangle^\bot\cap f^{-1}(\langle  {\mathbf 1} \rangle)= \bigcup_{\lambda \in \mathbb R}\langle {\mathbf 1} \rangle^\bot \cap f^{-1}(\lambda {\mathbf 1}),
    \label{eq:Ez_union}
\end{equation} 
where $\mathscr E_{\mathbf z}$ are the equilibria in the $\mathbf z$ coordinates and, 
\begin{equation}
    \langle {\mathbf 1} \rangle^\bot \cap f^{-1}(\lambda {\mathbf 1})=\left \{  {\mathbf z}\in \mathbb R^n\textrm{ : } \sum_{i=1}^nz_i =0 \textrm{ and } z_i\in \hat f^{-1}(\lambda)\right\}.
    \label{eq:Ez_one_elemente_union}
\end{equation}
Now if the sets $\langle  {\mathbf 1} \rangle^\bot \cap f^{-1}(\lambda {\mathbf 1})$  are non-empty for all $\lambda$ and they do not coincide, we expect to have a curve of equilibrium points parametrized by $\lambda$. Note that  $\langle  {\mathbf 1} \rangle^\bot \cap f^{-1}(\lambda {\mathbf 1})\not = \emptyset$ is equivalent to the existence of  $\alpha_1^\lambda,\dots,\alpha _n^\lambda $ roots of $\hat f -\lambda $ such that $\sum_{i=1}^n\alpha _i^\lambda=0$. The following polynomials are an important class of functions that satisfy this condition. 

\begin{proposition}
Consider \cref{principal} for the network $\mathcal C_n$ and $f_{ij}=\hat f$ for all $\{i,j\}\in \mathcal E$. Assume that $\hat f$ is an odd polynomial of degree $k\not =1$ with $k$ distinct real roots and $k|n$. Then, \cref{principal} has a continuum of equilibrium points. 
\label{prop:continuom_equilibria_polinomial}
\end{proposition}
\begin{proof}
Denote by $a_0,\dots,a_k$ the coefficients of $\hat f$. As $\hat f$ has $k$ distinct roots and all of them are simple, there exists $\epsilon>0$ such that for all $\lambda \in (-\epsilon, \epsilon)$, $\hat f-\lambda$ has $k$ real roots which we denote by  $\beta_0^\lambda,\dots, \beta_{k-1}^\lambda$. Moreover, for each $i$, $\beta_i^\lambda$ takes a continuum of values parametrized by $\lambda$. As $\hat f$ is an odd function, its even coefficients are null, and in particular $a_{k-1}=0$. By the well known Vieta's formulas \cite[p. 99]{Algebra:vol1}  we have for all $\lambda$,
\[0=-\frac{a_{k-1}}{a_k}=\sum_{i=0}^{k-1}\beta_i^\lambda.\]
Now for $i\in \{1,\dots,n\}$ define $\alpha_i^\lambda= \beta_{i \hspace{1.5pt}\textrm{mod}\hspace{1.5pt} k}^\lambda$,  and note that   $\sum_{i=1}^n\alpha_i^\lambda=\frac{n}{k}(\sum_{i=0}^{k-1}\beta_i^\lambda)=0$ as $k|n$. So from \cref{eq:Ez_one_elemente_union} and \cref{eq:Ez_union}, we deduce that 
\[(\alpha_i^\lambda)_{i=1}^n\in \langle {\mathbf 1} \rangle^\bot \cap f^{-1}(\lambda {\mathbf 1})\subset \mathscr E_{\mathbf z}\]
for all $\lambda \in (-\epsilon, \epsilon)$. Hence, $(\alpha_i^\lambda)_{i=1}^n$ gives a continuum of equilibria parametrized by $\lambda$.
\end{proof}

\subsection{Complete graphs or cliques}
\label{sec:complete_graph}
In this section we will deal with the complete graph of $n$ nodes $\mathcal K_n$. We will limit our study to the subset of equilibrium points $\tilde {\mathscr E}_{\mathbf x}$ with $f_{ij}=\hat f $ for all $\{i,j\}\in \mathcal E$. Moreover, it will be more convenient to think of the state space as $\mathbb R^n/\langle\mathbf 1\rangle $ rather than $\langle \mathbf 1\rangle ^\bot$. 

First, assume that $\hat f $ only has 3 distinct roots, which we denote by $0$, $\pm \alpha$ and let $\mathbf x^\star \in \tilde {\mathscr E}_{\mathbf x}$. Then, by definition of $\tilde {\mathscr E}_{\mathbf x}$ on a complete graph, $x_i^\star -x^\star_j\in  \{0,\pm\alpha\}$ for all $i,j$.  With elementary arguments (see \cref{supp:sec:additive_open}) one can show that for an appropriate representative and ordering of the nodes, 
 \begin{equation}
   \mathbf  x^\star =(0,\dots,0,\alpha,\dots,\alpha)^\top,
   \label{(0...0.alpha...alpha)}
\end{equation}
with $1$ up to $n$ null entries, depending on the equilibrium point.
In fact, the equilibria in $\tilde {\mathscr E}_{\mathbf x}$ will be of the form \cref{(0...0.alpha...alpha)} even if $\hat f$ is an arbitrary odd function,  as long as  $\hat f ^{-1}(0)\cap (0,\infty)$ is ``additive open'', i.e. $a+b\not =c$ for all $a,b,c\in \hat f ^{-1}(0)\cap (0,\infty)$, see \cref{supp:sec:additive_open}. We characterise the stability of these equilibrium points. 


\begin{proposition}
Let $\alpha $ be a root of $\hat f $ and consider \cref{principal} in $\mathcal K_n$ with $n\geq 3$. Assume that $\hat f'(0)> 0$, $\hat f '(\alpha)<0$ and denote
\[a=\frac{\hat f '(0)}{\hat f '(0)+|\hat f '(\alpha)|}, \hspace{1cm} b=\frac{|\hat f '(\alpha)|}{\hat f '(0)+|\hat f '(\alpha)|}.\]
Then, the stability of $\mathbf x^\star$ as in   \cref{(0...0.alpha...alpha)} with $n_0$ zero entries is given by, 
\begin{align*} 
    {n_0}/{n}&\in \left (a,b\right )\Longrightarrow \textrm{ $\mathbf x^\star$ is stable} ,\\ 
    {n_0}/{n}&\not \in \left [a,b\right ]\Longrightarrow \textrm{ $\mathbf x^\star$ is unstable}.
\end{align*}
\label{prop:kn_stability}
\end{proposition}

    Our proof follows from the direct computation of the Jacobian eigenvalues for a given number of zero entries in $\mathbf{x}^\star$, see  \cref{app:kn}. Moreover, one gets the same instability condition by applying  \cref{thm:unstability_r^-_ij r^+_ij} in this context (see  \cref{supp:sec:instability_kn}), which shows how powerful this theorem can be. 

\begin{remark}
With the notation of the previous proposition, if $\hat f '(0),\hat f '(\alpha)>0$ (resp. $\hat f '(0),\hat f '(\alpha)<0$) we get that $\mathbf x^\star$ is unstable (resp. stable), by \cref{prop:G-connected}.  If $\hat f '(0)<0$ and $\hat f '(\alpha)> 0$,  then $\mathbf x^\star$ is stable if $n_0\in \{ 0,n\}$ and unstable otherwise, again by \cref{prop:G-connected}. 

\end{remark}


\subsection{Example of the study of a tree of motifs}
\label{supp:sec:example}
\begin{figure}[htbp]
    \centering
        {\includegraphics[clip, trim=0cm 5.8cm 0cm 6.3cm, width=0.99\textwidth]{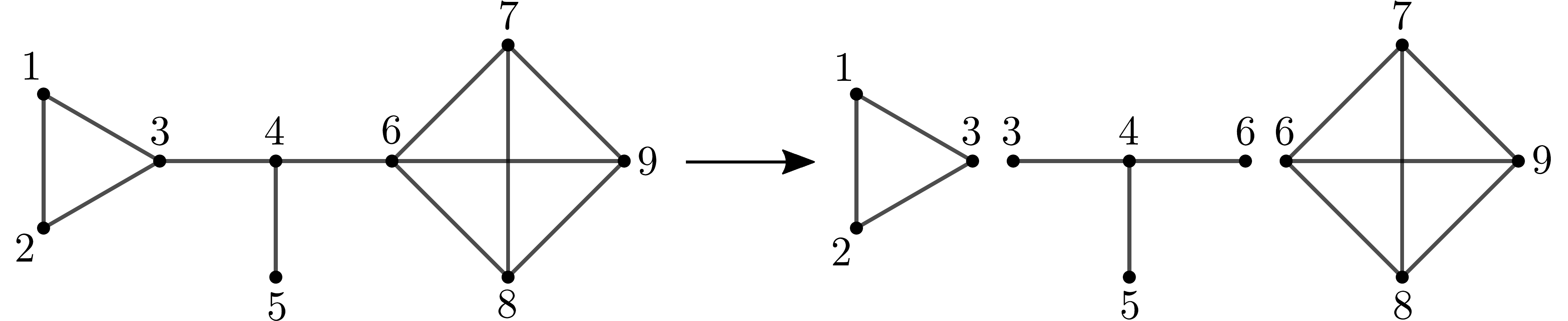}}
    \caption{On the left hand side the network $\mathcal G$ is depicted with its nodes labeled. On the right-hand side we depict the relevant subnetworks to study by  \cref{section:union_graphs}.}
    \label{supp:fig:geo_graph}
\end{figure}

Consider the system \cref{principal} with $f_{ij}(y)=\hat f(y)=y-y^3$ for all $i,j$, in the network $\mathcal G$ depicted in \cref{supp:fig:geo_graph}. By the main result of  \cref{section:union_graphs} we can find its equilibria and their stability by separately studying the subnetworks showed on the right-hand side of \cref{supp:fig:geo_graph}, i.e. a cycle graph $\mathcal C_3$, a star graph with 3 edges and a complete graph $\mathcal K_4$. 


\textbf{Cycle graph $\mathcal C_3$:} The results from  \cref{fig:circle} with $\lambda =1$ claim that, in $\langle \mathbf 1\rangle ^\bot$ this system has an unstable equilibrium at $\mathbf 0$ and a circle of stable ones, which can be parametrized by,
\[(x_1^\star,x_2^\star,x_3^\star)^\top=\frac{\sqrt{2}}{\sqrt{3}}\left( \frac{1}{\sqrt{2}}\cos(\theta)(-1,1,0)+\frac{1}{\sqrt{6}}\sin (\theta)(1,1,-2) \right)^\top, \]
for $\theta \in [0,2\pi)$. In the edge space the set of stable equilibria is given by,
\[\mathcal P_1^-=\left \{\frac{1}{\sqrt{3}}\left ( 2\cos(\theta), \cos(\theta)-\sqrt{3}\sin(\theta),-\cos (\theta)-\sqrt{3}\sin (\theta)\right )^\top\hspace{2pt}:\hspace{2pt}\theta\in [0,2\pi)\right \},\]
and the unstable ones are $\mathcal P_1^+=\{(0,0,0)^\top\}$.

\textbf{Star graph with 3 edges:} Note that the roots of $\hat f(y) = y - y ^3$ are $\{0,\pm 1\}$, and $\hat f'(0)=1$, $\hat f'(\pm 1)=-2$. So by \cref{prop:tree_stability} the stable equilibria in the edge space are given by, 
\[\mathcal P_2^- =\{\pm 1\}^3 :=\{\pm1\}\times\{\pm1\}\times\{\pm1\},\]
and the unstable ones by,
\[\mathcal P_2^+ =\left \{\left (y_{\{3,4\}},y_{\{4,5\}},y_{\{4,6\}}\right )^\top\in \{0,\pm 1\}^3 \hspace{2pt}:\hspace{2pt} y_{\{3,4\}}y_{\{4,5\}}y_{\{4,6\}}=0 \right \}.\]

\textbf{Complete graph $\mathcal K_4$:} Using the results from  \cref{sec:complete_graph}, and imposing the condition of orthogonality to $\mathbf 1$ we find that $\mathbf x^\star =(x_6^\star,x_7^\star,x_8^\star, x^\star_9)^\top\in \tilde{ \mathscr  E}_{\mathbf x}$ if and only if, 
\begin{equation}
    \{x_6^\star,x_7^\star,x_8^\star,x_9^\star\}\in\left \{\left\{0,0,0,0\right\},
    \left\{\tfrac{-1}{4},\tfrac{-1}{4},\tfrac{-1}{4},\tfrac{3}{4}\right\},
    \left\{\tfrac{1}{4},\tfrac{1}{4},\tfrac{1}{4},\tfrac{-3}{4}\right\},
    \left\{\tfrac{1}{2},\tfrac{1}{2},\tfrac{-1}{2},\tfrac{-1}{2}\right\}\right \},
    \label{supp:eq:eq_k4}
\end{equation} 
where we use the brackets $\scalebox{0.87}{\{\}}$ to denote multisets. 
To find the stability of these equilibria recall the notation of \cref{prop:kn_stability} and notice that in our case,
\[a=\frac{\hat f '(0)}{\hat f '(0)+|\hat f '(1)|}=\frac{1}{3}, \hspace{1.5cm}b=\frac{|\hat f '(1)|}{\hat f '(0)+|\hat f '(1)|}=\frac{2}{3}.\]
Now, note that the representative in the form of \cref{(0...0.alpha...alpha)} 
of the  first three types of equilibria in \cref{supp:eq:eq_k4}, have respectively 4, 3 and 1 null entries.  Hence, for these points $n_0/n\not \in [a,b]$ and thus, they are unstable. We denote by $\tilde {\mathcal P}_3^+$ the set of these points in edge coordinates. For the forth type of equilibria in \cref{supp:eq:eq_k4} the representative has two null entries, so $n_0/n=1/2\in (a,b)$ and thus, these points are stable. We denote by $\mathcal P_3^-$ the set of these points in the edge coordinates. 

Recall that for the complete graph there may be equilibria outside of  $\tilde {\mathscr E}_{\mathbf x}$. As we are working with polynomials over the integers, we can hope to find all equilibria using Gr\"obner Bases (see \cite{grobner}). Using \texttt{Singular}, we find that there is only one other type of equilibria,
\[ \{x_6^\star,x_7^\star,x_8^\star,x_9^\star\}=\left \{0,0,\tfrac{\sqrt{2}}{\sqrt{5}},\tfrac{-\sqrt{2}}{\sqrt{5}}\right \}.\]
As this type of equilibrium point is not in $\tilde {\mathscr E}_{\mathbf x}$, we can not use \cref{prop:kn_stability} to determine its stability. However, if we let $i$ and $j$ be the nodes with value 0, one finds that $r_{ij}^+=1$ and $r_{ij}^-=5$. So we have $r_{ij}^->r_{ij}^+$ and by  \cref{thm:unstability_r^-_ij r^+_ij} we conclude that these points are unstable. Denote by $\mathcal P_3^+$ the union of these points in edge coordinates with $\tilde {\mathcal P}_3^+$.

We are now prepared to tackle the whole network $\mathcal G$. For convenience, we denote,
\[\mathbf y =\left (\mathbf y_1, \mathbf y_2 ,\mathbf y_3 \right )^\top, \]
where, 
\[\mathbf y_1=\left(y_{\{1,2\}},y_{\{1,3\}},y_{\{2,3\}}\right )^\top, \hspace{1.5cm} \mathbf y_2= \left (y_{\{3,4\}},y_{\{4,5\}},y_{\{4,6\}}\right )^\top,\]
and,
\[\mathbf y_3=\left(y_{\{6,7\}},y_{\{6,8\}},y_{\{6,9\}},y_{\{7,8\}},y_{\{7,9\}},y_{\{8,9\}} \right)^\top.\]
Let $\mathcal P_i=\mathcal P_i^+\cup \mathcal P_i^-$ for $i=1,2,3$. Then, by the results from  \cref{section:union_graphs}, the stable equilibrium points in edge coordinates are, 
\[\mathcal P_-=\{(\mathbf y_1, \mathbf y_2 ,\mathbf y_3 )^\top \hspace{2pt}:\hspace{2pt} \mathbf y_i\in \mathcal P_i^- \},\]
and the unstable ones are,
\[\mathcal P_+=\{(\mathbf y_1, \mathbf y_2 ,\mathbf y_3 )^\top \hspace{2pt}:\hspace{2pt} \mathbf y_i\in \mathcal P_i \textrm{ and } \mathbf y\not \in \mathcal P_- \}.\]
 One can use the map $d$ to find the node coordinates of these states. Note that the node coordinates in the whole network will not coincide with the node coordinates found for the different subnetworks.

\section{Conclusion}
In this paper, we have presented an in-depth study of a general nonlinear model of consensus dynamics, showing a rich phenomenology depending on the structural properties of its underlying network. The dynamical model can be viewed as a gradient system that conserves the mean state. Furthermore, it does not produce complex dynamical structures such as periodic orbits, as  all orbits converge to equilibria.  Notably for a high-dimensional nonlinear system with arbitrary topology, we find a compact expression for the equilibrium points that highlights how these equilibria are determined by an interplay between the coupling function and the underlying network. Moreover, in line with previous results in the literature, we find that the stability of certain equilibria depends on effective resistances in the network, where the influence of general coupling functions can be taken into account using appropriate link weights.

Our analysis provides insight on simple networks like trees, cycles and cliques, but also for any combination of them in a tree of motifs, as we show that knowledge of the equilibria and their stability for individual motifs is in that case sufficient to characterise the whole-network dynamics. Although these conditions may seem  restrictive, the resulting structures may be seen as a generalisation of trees to the case of higher-order networks \cite{lambiotte2019networks}. Moreover, locally tree-like structures of cliques are expected to appear in \emph{projections of bipartite graphs} \cite{bipartite_projection}, which include co-author networks, and  when modeling  \emph{pervasive overlap} \cite{nature_figure,bodo2016sis} in social networks.


As a perspective for future research, we believe that a more careful deduction in the line of  \cref{thm:resistencia:one_edge} and   \cref{lem:instabiliy} could lead to better stability conditions for arbitrary equilibrium points, compared to those now presented in \cref{thm:resistencia:mult_edge}. Also in the spirit of \cref{thm:resistencia:one_edge}, one could specialize the Schur complement reduction to particular configurations of positive edges, leading to tight stability conditions. As a third extension, it might be possible to obtain stability results for large dense networks (which are ``approximately'' complete) based on the exact Jacobian eigenvalues for complete graphs (as in \cref{app:kn}) and invoking perturbation arguments on the Jacobian \cite{LotkerZvi2007Noda}. A similar approximate setting where some of our exact results might be used as a starting point is the study of networks whose local structure may be approximated by a tree of motifs, generalising standard approximations based on a locally tree-like structure \cite{melnik2011unreasonable}.  
Finally, while this work is theoretical in nature we believe that our developed insights can be used for the application and specialisation of system \cref{principal}  in a practical context. 

\appendix

\section{Proof of  \texorpdfstring{\cref{prop:G-connected}}{ref{prop:G-connected}}}
\label{app:connected}
For the stability condition, note that we have $J=-L^-$. Now as mentioned below \cref{eq:laplacian_matrix_prod}, $L^-$ is positive semi-definite and, as $\mathcal G^-$ is connected, 0 has multiplicity one in $\mathbb R^n$. Thus, in the state space $\langle \mathbf 1\rangle^\bot$ all eigenvalues of $-L^-$ are negative.

If $\mathcal G^-$ is disconnected, consider $\mathcal V_1$ a connected component and the vector $\mathbf x= \sum_{i\in \mathcal V_1} \mathbf e_i$. Then, follow the proof of \cref{prop:cut_set_unstability} to deduce that $\lambda_{\max}\geq 0$.   If there is an edge between $\mathcal V_1$ and $(\mathcal V\setminus \mathcal V_1)$ in $\mathcal G^+$, we may directly apply \cref{prop:cut_set_unstability} to the associated cut-set. 

\section{Further Schur complement reduction}
\label{app:Schur}
\begin{figure}[htbp]
    \centering
        \includegraphics[clip, trim=3.3cm 0.95cm 6.2cm 0.88cm, width=0.95\textwidth]{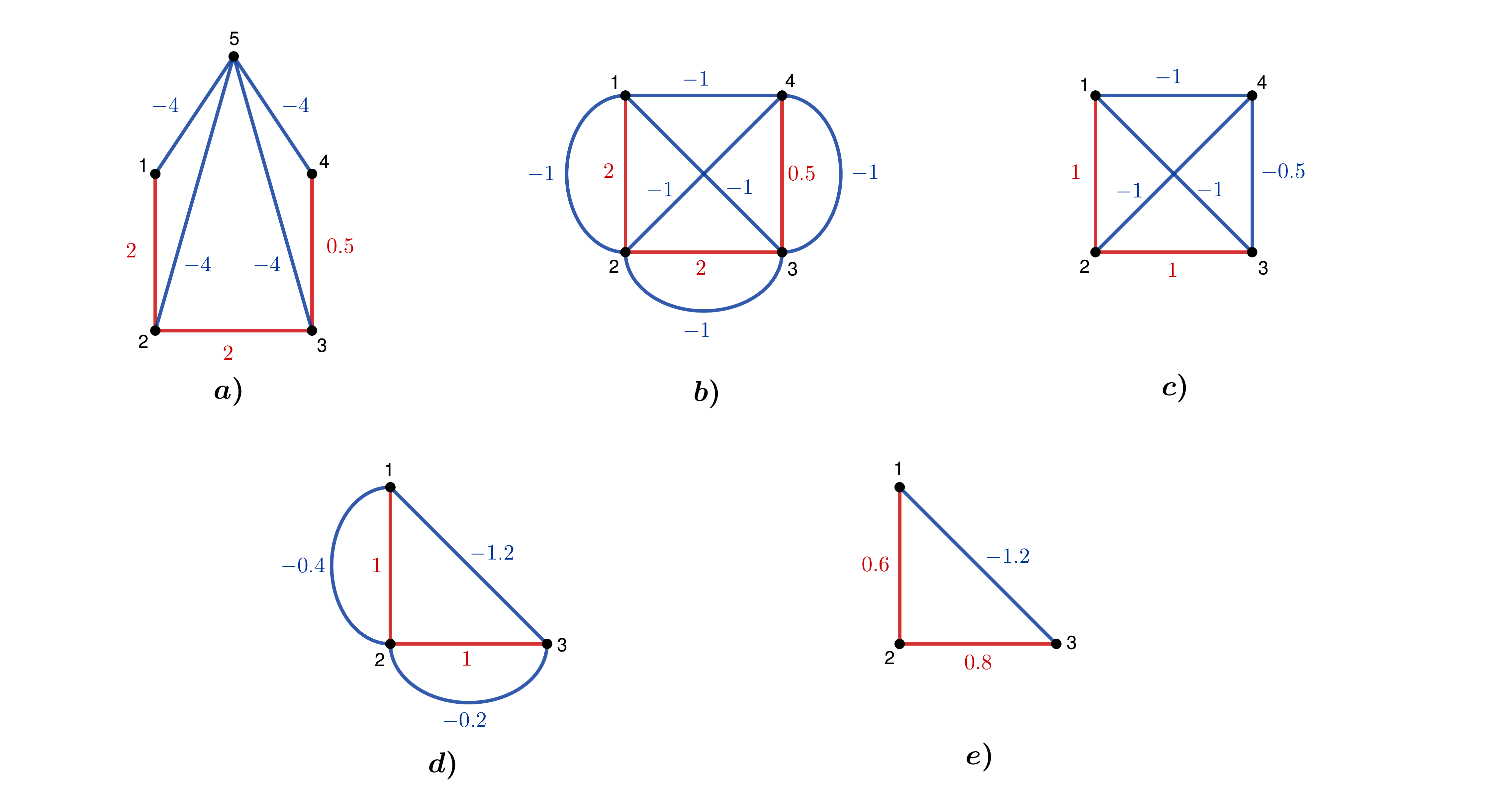}
    \caption{Networks obtained by the recursive application of Schur complement. We depict edges from the corresponding $L^+$ in \emph{red}. Edges from $L^-$ are depicted in \emph{blue} and with a negative sign added to their weights.}
    \label{supp:fig:schur_comp}
\end{figure}

In this section we explain how to apply the Schur complement reduction recursively through an example. Consider a Jacobian matrix $J=L^+-L^-$, where the weights of $L^+$ and $L^-$ are depicted in \cref{supp:fig:schur_comp}\textcolor{My_green10}{.a}. Note that there is a unique node not incident to any edge in $\mathcal G^+$, i.e $\mathcal N= \{5\}$. By \cref{thm:schur_comp}, the linear stability of $J$ is determined by $J/\mathcal N= L^+_{\mathcal N^c,\mathcal N^c}- L^-/\mathcal N$, where the weights of these Laplacians are depicted in \cref{supp:fig:schur_comp}\textcolor{My_green10}{.b}. When we have multiple edges between a pair of nodes in \cref{supp:fig:schur_comp}\textcolor{My_green10}{.b} we can simply add the (signed) weights to get \cref{supp:fig:schur_comp}\textcolor{My_green10}{.c}. This corresponds to finding another decomposition  $J/\mathcal N= \tilde L^+- \tilde L^-$ with Laplacians which do not have common edges. We observe that in \cref{supp:fig:schur_comp}\textcolor{My_green10}{.c}, the node $4$ is not incident to any edge in   $\tilde L^+$. Thus, we can apply our result again with $\tilde {\mathcal N}=\{4\}$, to get \cref{supp:fig:schur_comp}\textcolor{My_green10}{.d} and  by adding edges we get \cref{supp:fig:schur_comp}\textcolor{My_green10}{.e}. In summary, we have shown that the linear stability of $J$ is given by the linear stability of $(J/\mathcal N)/\tilde {\mathcal N}$ depicted in \cref{supp:fig:schur_comp}\textcolor{My_green10}{.e}, which is unstable by \cref{prop:G-connected}. 

Note that in this example we have considered $(J/\mathcal N)/ \tilde {\mathcal N}$. The quotient formula of Schur complement  \cite[Equation 6.0.26]{Schur_complement} states that if $\mathcal N$ and $\tilde {\mathcal N}$ are disjoint sets of nodes then 
\[(J/\mathcal N)/ \tilde {\mathcal N}= J/(\mathcal N\cup \tilde {\mathcal N}), \]
when the left-hand side is well defined. 
Given a Jacobian of our system it would be interesting to be able to determine which nodes we will be able to remove applying the Schur complement process recursively. Then, by the quotient formula we could reduce all of them at once. One can sideline this problem by taking the opposite approach. That is, in each iteration of the Schur reduction, only reduce one node until no node is isolated in the corresponding $\mathcal G^+$. This approach has the advantage that each iteration only requires inverting a scalar instead of a matrix, see \cref{eq:schur_def}.

\section{Proof of \texorpdfstring{\cref{thm:resistencia:mult_edge}}{ref{thm:resistencia:mult:edge}}}
\label{app:resistance_thm} 
For the instability note that if $\mathcal G^-$ is disconnected, the result follows from \cref{prop:G-connected}. Otherwise, we can write $J = L^+_{ij} + \tilde L^+ -L^-$, where $L^+_{ij}$ is the Laplacian with only edge $\{i,j\}$ with weight coming from $\mathcal G^+$ and $\tilde L^+$ is the Laplacian of the rest of edges from $\mathcal G^+$. Now as $ \tilde L^+$ is positive semi-definite, it is enough to show that $L^+_{ij} -L^-$ is linearly unstable, which follows from \cref{rem:one_edge}.

For the stability result, as $\sum_{\{i,j\}\in \mathcal G^+ }f'_{ij}(x^\star_i-x^\star_j)r_{ij}^-<1$  we have $r_{ij}^-<\infty$ for all edges in $\mathcal G^+$, and thus $\mathcal G^-$ is connected.
Moreover, we can take weights $\alpha_{ij}>f'_{ij}(x^\star_i-x^\star_j)r_{ij}^-$ such that $\sum_{\{i,j\}\in \mathcal G^+ } \alpha_{ij}=1$. Then, using the notation $L^+_{ij}$ from above we have,
\[J= L^+-L^-=\sum_{\{i,j\}\in \mathcal G^+ }L_{ij}^+-\alpha_{ij} L^- = \sum_{\{i,j\}\in \mathcal G^+ }\alpha_{ij}(\alpha_{ij}^{-1}L_{ij}^+- L^-).
\]
Now apply \cref{rem:one_edge} to each $\alpha_{ij}^{-1}L_{ij}^+- L^-$ to conclude that they are linearly stable matrices, thus so is $J$ and    $\mathbf x^\star$ is stable.

\section{Proof of \texorpdfstring{\cref{lem:instabiliy}}{ref{lem:instabiliy}}}
\label{app:lem:unstable}
First, assume that $A$ and $B$ are non-singular, and let $B^{1/2}$ be the positive definite matrix such that $B=B^{1/2}B^{1/2}$. Then, 
\begin{align*}
\max_{\mathbf x}\frac{\mathbf x^\top A\mathbf x}{\mathbf x^\top B\mathbf x}&=\max_{\mathbf x}\frac{\mathbf x^\top  A\mathbf x}{(B^{1/2}\mathbf x)^\top (B^{1/2}\mathbf x)}= \max_{\mathbf y}\frac{\mathbf y^\top  B^{-1/2} A B^{-1/2}\mathbf y}{\mathbf y^\top \mathbf y} \\
&=\max_{||\mathbf y||=1}{\mathbf y^\top  B^{-1/2} A B^{-1/2}\mathbf y}=\mu_{\max},
\end{align*}
where $\mathbf y = B^{1/2}\mathbf x$, $\mu_{\max}$ is the largest eigenvalue of $B^{-1/2} A B^{-1/2}$ and in the last equality we use the Courant minmax principle \cite{LotkerZvi2007Noda}. Now, note that $\mu_{\max}^{-1}=\tilde \mu_{\min}$ where $\tilde \mu_{\min}$ is the minimum eigenvalue of $B^{1/2} A^{-1} B^{1/2}$, so 
\begin{align*}
\max_{\mathbf x}\frac{\mathbf x^\top  A\mathbf x}{\mathbf x^\top B\mathbf x}&=\mu_{\max}=\frac{1}{\mu_{\max}^{-1}}=\frac{1}{\tilde \mu_{\min}}=\frac{1}{\min_{||\mathbf y||=1} \mathbf y^\top B^{1/2} A^{-1} B^{1/2}\mathbf y} \\
&=\max_{||\mathbf y||=1}\frac{1}{ \mathbf y^\top B^{1/2} A^{-1} B^{1/2}\mathbf y}=\max_{\mathbf y}\frac{\mathbf y^\top \mathbf y}{ \mathbf y^\top B^{1/2} A^{-1} B^{1/2}\mathbf y}= \max_{\mathbf z} \frac{\mathbf z^\top  B^{-1}\mathbf z}{\mathbf z^\top A^{-1}\mathbf z}.
\end{align*}


For the general case, let $V=(\ker A + \ker B)^\bot $ and note that $A$ and $B$ are positive definite in this space. Recall from the preliminaries that  $(A|_V)^{-1}=A^\dag|_V$ (resp. for B), so from the argument above we get,
\begin{equation*}
 \max_{\mathbf x\in V}\frac{\mathbf x^\top  A\mathbf x}{\mathbf x^\top B\mathbf x}=\max_{\mathbf x\in V}\frac{\mathbf x^\top  B^\dag \mathbf x}{\mathbf x^\top A^\dag \mathbf x}.
\end{equation*}
Now as $A,B$ are positive semi-definite and  $\ker A=\ker A^\dag$ (resp. for $B$), the previous expression is equivalent to, 
\begin{equation*}
 \max_{\mathbf x\bot \ker B}\frac{\mathbf x^\top  A\mathbf x}{\mathbf x^\top B\mathbf x}=\max_{\mathbf x\bot \ker A}\frac{\mathbf x^\top  B^\dag \mathbf x}{\mathbf x^\top A^\dag \mathbf x}.
\end{equation*}

\section{Proof of \texorpdfstring{\cref{union_graph:lema}}{ref{union:graph:lema}}}
\label{app:union_graphs}
Reordering the nodes we can assume that the deleted row in $\tilde d$ is the first one.  Denote by $\lambda_1\leq \dots\leq  \lambda_n$ the eigenvalues of $dMd^\top $ and by $\tilde \lambda_1\leq \dots\leq \tilde \lambda_{n-1}$ the ones of $\tilde dM\tilde d^\top$. As $dMd^\top $ is a symmetric matrix and $\tilde dM\tilde d^\top$ is the same matrix with the first row and column deleted we have, 
\begin{equation}
    dM d^\top=\begin{pmatrix}
a &\mathbf b^\top \\
\mathbf b&\tilde dM\tilde d^\top
\end{pmatrix}.
\label{eq:lem:dMd}
\end{equation} 
Thus we can apply Cauchy's Interlace theorem \cite{Cauchy_interlace} and we get, 
\begin{equation}
    \lambda _1\leq \tilde \lambda_1\leq \lambda _2\leq  \dots \leq \lambda _{n-1}\leq \tilde \lambda_{n-1}\leq \lambda _n.
    \label{eq:interlace:lem}
\end{equation}
Recall that $d^\top \mathbf 1=0$, so the multiplicity of the eigenvalue  0 for $dMd^\top$, which we denote by $k$, is at least one. By \cref{eq:interlace:lem} it is clear that the multiplicity of $0$ for $\tilde dM\tilde d^\top$, which we denote by $\tilde k$, is $k-1$, $k$ or $k+1$. If $\tilde k= k-1$, using \cref{eq:interlace:lem} it is clear that the number of positive (resp. negative) eigenvalues of each matrix coincide. We now show by contradiction that $\tilde k\not = k , k+1$. 

First note that by \cref{eq:lem:dMd} for all  $\mathbf y\in \ker \tilde  dM\tilde d^\top \cap \langle \mathbf b\rangle ^\bot $, $(0,\mathbf y^\top )^\top \in \ker dMd^\top $. Moreover, $\mathbf 1\in \ker dMd^\top $, so
\begin{equation}
    k= \textrm{dim}(\ker dMd^\top)\geq \textrm{dim}\left ( \ker \tilde  dM\tilde d^\top \cap \langle \mathbf b\rangle ^\bot \right )+1.
    \label{eq:lema:dim}
\end{equation}
Now, if $\tilde k= k+1$, then $\ker \tilde  dM\tilde d^\top \cap \langle \mathbf b\rangle ^\bot  $ has at least dimension $k$ so we get a contradiction with \cref{eq:lema:dim}. If $\tilde k= k$, the multiplicity of 0 in both matrices coincide, thus we can apply an extension of Cauchy Interlace Theorem (see \cite[Theorem 2]{Cauchy_interlace}), which claims that $\langle\mathbf b\rangle\bot \ker \tilde dM \tilde d^\top$. Then, $\ker \tilde  dM\tilde d^\top \cap \langle \mathbf b\rangle ^\bot =\ker \tilde  dM\tilde d^\top $ and again we get a contradiction with \cref{eq:lema:dim}.

\section{Proof of \texorpdfstring{\cref{prop:kn_stability}}{ref{prop:kn:stability}}}
\label{app:kn}
Our arguments closely resembles the one given in \cite{karel2020arXiv} for a specific function $\hat f$. 

Note that $a,b\in (0,1)$, so if $n_0\in \{0,n\}$ then $n_0/n\not \in [a,b]$. In these cases, for all edges $x_i^\star-x_j^\star=0$ and as $f'(0)>0$, $\mathbf x^\star$ is unstable by \cref{prop:G-connected}. Otherwise, recall that  $J= L^+-L^-$ and it is easy to check that,
\[ L^+=f'(0) \begin{pmatrix}n_0P_{n_0} & \mathbf 0_{n_0,n-n_0}\\\mathbf 0_{n-n_0,n_0}&(n-n_0)P_{n-n_0} \end{pmatrix} ,\hspace{22pt} L^- =-f'(\alpha)\begin{pmatrix}(n-n_0)I_{n_0} & -\mathbf{1}_{n_0}\mathbf{1}_{n-n_0}^\top\\-\mathbf{1}_{n-n_0}\mathbf{1}_{n_0}^\top&n_0I_{n-n_0} \end{pmatrix},\]
 where $P_{k}= I_k-\mathbf{1}_k\mathbf{1}_k^\top/k$, i.e. the orthogonal projection onto $\langle \mathbf 1 _k\rangle ^\bot$. Notice that this two matrices commute, hence they can be simultaneously diagonalized. We find the following four types of eigenvectors of $J$. 

\emph{Type 1:} $\mathbf 1_n$ is an eigenvector of eigenvalue $0$, which does not affect the stability in the state space $\langle \mathbf 1\rangle ^\bot$. 

\emph{Type 2:} $(-\mathbf 1_{n_0}^\top/n_0,\mathbf 1_{n-n_0}^\top/(n-n_0))^\top$ is an eigenvector of eigenvalue $\hat f '(\alpha)n$.

\emph{Type 3:} Any vector of the form $\mathbf z = (\mathbf z_{n_0}^\top,\mathbf 0_{n-n_0})^\top$ such that $\mathbf{1}_{n_0}^\top\mathbf z_{n_0} =0$ is an eigenvector of eigenvalue $(\hat f '(0)-\hat f '(\alpha))n_0+\hat f '(\alpha)n$.

\emph{Type 4:} Any vector of the form $\mathbf z = (\mathbf 0_{n_0},\mathbf z_{n-n_0}^\top)^\top$ such that $ \mathbf{1}_{n-n_0}^\top\mathbf z_{n-n_0}=0$ is an eigenvector of eigenvalue $(\hat f '(\alpha )-\hat f '(0))n_0+\hat f '(0)n$.

As there are $n_0-1$ linear independent eigenvectors of type 3 and $n-n_0-1$ eigenvectors of  type 4, these are all eigenvectors of $J$.

As $\hat f '(\alpha)n<0$, the stability of the equilibrium points is given by the signs of the other eigenvalues. Now  assume that $n_0\not \in \{ 1,n-1\}$ so that there is at least one eigenvector of each type. Type 3 eigenvalue is positive if $b<n_0/n$ and negative if $n_0/n<b$. Similarly, type 4 eigenvalue is positive if $n_0/n<a$ and negative if $a<n_0/n$. Thus, using \cref{prop:stability_equilibria} we get the desired result. 

If $n_0=1$ then we do not have eigenvectors of Type 3, so $\mathbf x^\star$ is stable if $a<n_0/n$ and unstable if $n_0/n<a$. Using that $n_0/n=1/n<1/2$ and  that $[a,b]\not = \emptyset$ if and only if $a\leq 1/2$, it follows,
\[a<n_0/n \Longleftrightarrow n_0/n\in (a,b), \hspace{1cm}\textrm{and}\hspace{1cm}n_0/n<a \Longleftrightarrow n_0/n\not \in [a,b].\]
Similar arguments work for the case $n_0=n-1$.

%

\bibliographystyle{siamplain}
\bibliography{references}
\end{document}


\maketitle

\section{Potential dynamics}
\label{supp:sec:dynamics}
We start by introducing some basic concepts of dynamical systems (we follow standard notation, see for instance \cite{supp_book_Gerald_Teschl}). We consider ODEs defined by $\mathcal C^1$ maps, 
\[\dot{\mathbf x}=g(\mathbf x),\]
and thus, by Picard–Lindelöf theorem \cite[Theorem 2.5]{supp_book_Gerald_Teschl} it will have existence and uniqueness of solution for any given initial condition. Moreover, to simplify some definitions of qualitative features, we will always assume that orbits are defined for all times. Given an ODE, $\dot{\mathbf x}=\tilde g(\mathbf x)$, this can can be done by considering the equivalent system,
\[\dot{\mathbf x}=g(\mathbf x):=\frac{\tilde g(\mathbf x)}{||\tilde g(\mathbf x)||+1},\]
which solutions are reparametrizations of solutions of the original ODE. 

 Let $\mathbf x_0\in \mathbb R^n$ and denote by $\varphi(t,\mathbf x_0)$ the solution of $\dot{\mathbf x}=g(\mathbf x)$ at time $t$ such that $\varphi(0,\mathbf x_0)=\mathbf x_0$. We denote its forward orbit and its $\omega$-limit by 
 \[\gamma _{+}(\mathbf x_0)=\{\varphi(t,\mathbf x_0)\textrm{ : } t>0\},\hspace{1cm}\omega (\mathbf x_0)=\left \{\lim_{n \to \infty} \varphi(t_n, x) \textrm{ : } t_n\xrightarrow[n \to \infty]{} +\infty\right \}.\]
 We say that an equilibrium point $\mathbf x^\star $ is \emph{stable} if for any neighbourhood $U$ of $\mathbf x^\star $, there exists a neighbourhood $V\subset U$ of $\mathbf x^\star $ such that for all $\mathbf x_0\in V$, $\gamma_+(\mathbf x_0)\subset U$. We say that  $\mathbf x^\star $ is \emph{attractive} if there exists a neighbourhood $U$ of $\mathbf x^\star$ such that for all $\mathbf x_0\in U$, $\omega(\mathbf x_0)=\mathbf x^\star$.
 Let $\mathbf x^\star $ be an equilibrium point (or infinity), we denote its basin of attraction, i.e. the set of points that converge to it, by $W^{+}(\mathbf x^\star)$. We define   $\gamma_-(\mathbf x_0)$, $\alpha(\mathbf x_0)$ and $W^{-}(\mathbf x^\star)$ in the same manner but for the system $\dot {\mathbf x}= -g(\mathbf x)$.  We will need the following result regarding the basin of attraction. 
 
\begin{lemma}
Consider $\dot {\mathbf x}=g(\mathbf x) $  where $g\in \mathcal C^1$, and let $\mathbf x^\star $ be an unstable hyperbolic equilibrium point. Then, $W^{+}(\mathbf x^\star)$ has measure\footnote{Lebesgue measure.} zero. 
\label{lemma:grad_mesure0}
\end{lemma}
\begin{proof}
 By the stable manifold theorem \cite[Theorem 9.4]{supp_book_Gerald_Teschl}, there is a neighbourhood $U$ of $\mathbf x^\star$, such that the set
\[W^{+}_U(\mathbf x^\star)=\{\mathbf x\in W^{+}(\mathbf x^\star)\textrm{ : }\gamma_+(\mathbf x)\cup \{\mathbf x\}\subset U \},\]
is a $\mathcal C^1$ submanifold of $\mathbb R^n$ tangent to the stable manifold of the linearization of the system. Since $\mathbf x^\star $ is not stable, these submanifolds have at most dimension $n-1$ and thus $W^{+}_U(\mathbf x^\star)$ has measure zero (e.g. by Sard's theorem \cite{supp_ArthurSard1942Tmot}). 
Then, note that,
\[W^{+}(\mathbf x^\star)=\bigcup_{m=0}^\infty \left \{\varphi(-m,\mathbf x)\textrm{ : }\mathbf x\in W^{+}_U(\mathbf x^\star)\right \}=\bigcup_{m=0}^\infty \varphi(-m,W^{+}_U(\mathbf x^\star)),\]
 and that,  for all $t$ the map $\varphi (t,\cdot )$ sends null sets to null sets, since it is a $\mathcal C^1$ diffeomorphism. So,  $W^{+}(\mathbf x^\star)$ is the countable union of null measure sets and thus, it has measure zero. 
\end{proof}

Now we focus on gradient dynamics.  Let $V:\mathbb R^n\rightarrow \mathbb R$ be a $\mathcal C^2$ function.  We define its \emph{gradient dynamics} as, 
 \begin{equation}
     \dot {\mathbf x} = -\nabla V(\mathbf x),
     \label{eq:Intro:gradient}
 \end{equation}
so that $V$ is a potential function of this system.
 The derivative of $V$ along the trajectories is, 
 \begin{equation}\dot V(\mathbf x)=\langle \nabla V(\mathbf x), \dot {\mathbf x}\rangle = - ||\nabla V(\mathbf x)||^2\leq 0,\label{eq:derivative_grad}\end{equation}
and note that  $\dot V(\mathbf x)=0$ only in equilibrium points. Thus, $V(\mathbf x(t))$ is decreasing in solutions of  \cref{eq:Intro:gradient} (except in equilibrium points). 

\begin{proposition}
Assume that the set of equilibrium points $\mathscr E$ of  the system \cref{eq:Intro:gradient} is countable. Then, if $\gamma _+(\mathbf x_0)$ (resp. $\gamma _-(\mathbf x_0)$)  is bounded, $\omega(\mathbf x_0) $ (resp. $\alpha(\mathbf x_0)$) is an equilibrium point.
\label{prop:grad}
\end{proposition}
\begin{proof}
As $\gamma _+(\mathbf x_0)$ is bounded, it is contained in a compact set so by LaSalle's invariance principle \cite{supp_LaSalleJosephP1976Tsod,supp_book_Gerald_Teschl},
 \[\omega(\mathbf x_0)\subset \{\mathbf x\in \mathbb R^n\textrm{ : }\dot V(\mathbf x)=0\}=\mathscr E.\]
Now  $\gamma _+(\mathbf x_0)$ is contained in a compact set so $\omega (\mathbf x_0)$ is non-empty and connected (see \cite[Lemma 6.6]{supp_book_Gerald_Teschl}). As $\mathscr E$ is countable the only non-empty connected subsets are singletons. 
\end{proof}

To be able to apply the previous result we need $\gamma _+(\mathbf x_0)$ to be bounded. One sufficient condition is for $V$ to be  \emph{radially unbounded}, i.e. $V(\mathbf x)\rightarrow +\infty$ when $||\mathbf x||\rightarrow +\infty$, as the potential is decreasing in orbits. Now we are prepared to show that, under mild conditions, any forward orbit of the dynamics \cref{principal} converges to a single equilibrium point.

\begin{proposition}
Let $\mathscr E_{\mathbf x}$ be the equilibrium points of \cref{principal} in the state space $\langle\mathbf 1 \rangle^\bot$. Assume that $\mathscr E_{\mathbf x}$ is countable and one of the following conditions holds:\emph{
\begin{enumerate}[label=(\alph*)]
    \item \emph{For all $\{i,j\}\in \mathcal E$, $-F_{ij}$ is radially unbounded.}
    \item \emph{For all $\{i,j\}\in \mathcal E$, $\limsup_{y\to+\infty}yf_{ij}(y)<0$.}
    \item \emph{All orbits are bounded,  i.e. contained in a compact set.}
    \end{enumerate}}
    Then, 
     \[\langle \mathbf 1\rangle ^\bot= \bigcup_{\mathbf x^\star\in  \mathscr E_{\mathbf x}}W^+(\mathbf x^\star)=W^-(\infty)\cup \bigcup_{\mathbf x^\star\in  \mathscr E_{\mathbf x}}W^-(\mathbf x^\star),\]
     where if \emph{(c)} holds, we also have $W^-(\infty)=\emptyset$.
     \label{prop:descomp_state_W}
 \end{proposition}
 \begin{proof}
By \cref{prop:grad} it is enough to show that forward orbits are bounded, and backwards unbounded ones tend to infinity. So condition (c) is clearly sufficient. 
 Condition (a) is also sufficient, as then   $V|_{\langle \mathbf 1 \rangle^\bot}$ defined by \cref{eq:def_grad} is radially unbounded. Finally, condition (b) is a particular case of (a), as we proceed to show.

Assume that (b) holds, and consider $g(y)=\sup_{z\geq y}zf_{ij}(z)$, which is a non-increasing function. Then, let 
\[a=\lim_{y\rightarrow +\infty}g(y),\]
and notice that by our assumption $a<0$. Thus, there exists $M>0$ such that $g(y)<\frac{a}{2}$ for all $y\geq M$. In particular, $f_{ij}(y)<\frac{a}{2} \frac{1}{y}$ for all $y\geq M$ and hence for certain constant $c$,
\[\lim_{y\rightarrow + \infty}-F_{ij}(y)=c+\int_M^\infty -f_{ij}(y)>c+ \int_M^\infty  - \frac{a}{2}\frac{1}{y} = +\infty,\]
as $- \frac{a}{2}>0$. Finally, $f_{ij}$ is odd so $F_{ij}$ is even and thus, $\lim_{|y|\rightarrow \infty} -F_{ij}(y)=+\infty$. 
\end{proof}
 
 


 Note that the coupled phase oscillators model with constants $w_i=0$, trivially satisfies (c) when studied in $\mathbb S^1\times \dots \times \mathbb S^1$.
 Another important case is given when  $f_{ij}=\hat f$ for all edges and   $\hat f$ has finitely many roots. In this case (and changing sense of time if necessary), the system satisfies condition (b), as long as $\hat f$ does not tend to 0 at $\infty$.

\section{Dimension of the set of equilibrium points}
\label{supp:sec:dim_equilibria}
We first show that, 
 \begin{equation}
     \textrm{dim}(\ker d) = m-n+1, \hspace{1cm}\textrm{and}\hspace{1cm}\textrm{dim}\left ((\ker d)^\bot \right )= n-1.
\label{eq:dim_cut}
 \end{equation}
 where  $n=|\mathcal G|$ and $m=|\mathcal E|$. To do so we find a basis of each space through a spanning tree\footnote{That is, a subnetwork of $\mathcal G$ with all of its nodes and no cycles.} $\mathcal T$  of $\mathcal G$, which also gives us some insight on how to explicitly compute these spaces. Each edge in $\mathcal G\setminus \mathcal T$ determines a unique cycle with this edge and edges of $\mathcal T$. These cycles form a basis of the cycle space \cite{supp_book_Alg_graph}. Similarly, each edge in $\mathcal T$ determines a unique cut-set with this edge and edges in $\mathcal G\setminus\mathcal  T$. These cut-sets form a basis of the cut space \cite{supp_book_Alg_graph}. In particular, we get  \Cref{eq:dim_cut}. 

Assume now that for all edges, $f_{ij}$ has null derivative in a finite set of points. Then, $f_{ij}$ has finite changes of monotonicity, so  the set $f_{ij}^{-1}(y)$ is finite for any $y$. Thus,  the set $f^{-1}(\mathbf y)$ is also finite for any point $\mathbf y$. Informally, we can expect then that $f^{-1}(\ker d)$ is  $m-n+1$ dimensional, as $\dim (\ker d)=m-n+1$. Finally,  $\mathscr E_{\mathbf y}= (\ker d)^\bot \cap f^{-1}(\ker d  )$, so we expect the set of equilibria in $\langle \mathbf 1\rangle ^\bot$ to be at most $m-n+1$ dimensional. In \cref{sec:cycle_graph} we will show that for the cycle graph this bound may be reached and we get a continuum of equilibrium points in the node space, $\langle\mathbf 1\rangle^\bot$.

\section{Resistors in series and in parallel}
\label{supp:sec:parallel}
 The effective resistance as defined  in \cref{def:ef.res.} is precisely the effective resistance in the electrical circuit sense, where the circuit is given by the network and the resistors by the inverse of the edge weights. Thus, this definition satisfies the well known properties of resistors in series and in parallel, see \cite{supp_DorflerFlorian2018ENaA}. That is, if we have a series of nodes $i_1,\dots, i_k$ such that all paths from $i_1$ to $i_k$ go through $i_1,\dots, i_k$ in this order, then 
\begin{equation}
    r_{i_1i_k}= \sum_{j=1}^{k-1} r_{i_ji_{j+1}}.
    \label{supp:eq:series_res}
\end{equation}
For parallel resistors, if we have nodes $i,j$ and a partition  $\mathcal E=\mathcal E_1\cup \mathcal E_2$ such that all paths from $i$ to  $j$ are entirely contained in $\mathcal E_1$ or in $\mathcal E_2$ then, 
\begin{equation}
  \frac{1}{r_{ij}}=\frac{1}{r_{ij}^{\mathcal E_1}} +\frac{1}{r_{ij}^{\mathcal E_2}} ,
  \label{supp:eq:parallel_res}
\end{equation}
where $r_{ij}^{\mathcal E_1}$ is the effective resistance in the network with only the edges $\mathcal E_1$ (resp. for $\mathcal E_2$). 

The prototypical examples of resistances in series and parallel are depicted in the following diagrams. 
\begin{center}
\begin{circuitikz} \draw
(0,0) node[label={[font=\footnotesize]above:1}] {}to[resistor,*-*,l=$r_{1,2}$] (3,0) node[label={[font=\footnotesize]above:2}] {}
  to[resistor,*-*,l=$r_{2,3}$] (6,0)  node[label={[font=\footnotesize]above:3}] {}
  
  (9.5,-1.5) node[label={[font=\footnotesize]below:1}] {} to[short,*-] (8,-1.5) to[resistor,l=$r_{1,2}^{\mathcal E_1}$]  (8,1.5)
  to[short,-*] (9.5,1.5)node[label={[font=\footnotesize]above:2}] {} to (11,1.5) to[resistor,l=$r_{1,2}^{\mathcal E_2}$]  (11,-1.5)
  to (9.5,-1.5)
;
\end{circuitikz}
\end{center}
Note that the underlying network on the second diagram is a multigraph and not a graph. This is not relevant as effective resistance can be defined similarly for these.

\section{Using \texorpdfstring{\cref{lem:instabiliy}}{ref{lem:instabiliy}} for stability condition}
\label{supp:sec:stability B^T/A^T} 
We assume that $\mathcal G^-$ is connected, as by \cref{prop:G-connected} it is a necessary condition to have linear stability. Then, by \ref{prop:minimax} let $\mathbf z\in \langle \mathbf 1\rangle ^\bot$ be unitary such that,
\[\lambda_{\max} = \mathbf z^\top L^+\mathbf z-\mathbf z^\top L^-\mathbf z= \mathbf z^\top L^-\mathbf z\left (\frac{\mathbf z^\top L^+\mathbf z}{\mathbf z^\top L^-\mathbf z}-1\right ), \]
where we have used that,  $L^-$ is positive definite in our space $\langle \mathbf 1\rangle^\bot $ (as $\mathcal G^-$ is connected), so $\mathbf z^\top L^-\mathbf z>0$. Thus\footnote{As in the main article, all maximums are taken over non-null vectors.}, 
\[\textrm{sign}(\lambda_{\max })=\textrm{sign}\hspace{-2pt}\left (\frac{\mathbf z^\top L^+\mathbf z}{\mathbf z^\top L^-\mathbf z}-1\right )\hspace{-2pt}=\textrm{sign}\hspace{-2pt}\left (\max_{\mathbf x\bot \mathbf 1} \frac{\mathbf x^\top L^+\mathbf x}{\mathbf x^\top L^-\mathbf x}-1\right )\hspace{-2pt}=\textrm{sign}\hspace{-2pt}\left (\max_{\mathbf x\bot \ker {L^+}} \frac{\mathbf x^\top (L^-)^\dag\mathbf x}{\mathbf x^\top (L^+)^\dag\mathbf x}-1\right ),\]
where in the last equality we use  \cref{lem:instabiliy}.  We believe this last expression can lead to stability criteria. For instance, if $i$ is an isolated node in $\mathcal G^+$ and $\mathbf x\bot \ker L^+$, then $[\mathbf x]_i=0$ (as  $\ker L^+$ is generated by the connected components of $\mathcal G^+$). This gives rise to a different proof of the stability reduction \cref{thm:schur_comp}. Moreover,  $\mathbf x\bot \ker L^+$ imposes more conditions on $\mathbf x$, coming from the connected components of $\mathcal G^+$ with more than one node, which could potentially be used to reduce our problem even further.  

Another approach is to use the fact that the pseudoinverse of a Laplacian is closely related to the effective resistance, even more than the  Laplacian itself.
This might lead to stability criteria with effective resistance, as we did for instability in  \cref{thm:unstability_r^-_ij r^+_ij}.

\section{Equilibrium points of the coalescence of two networks}
\label{supp:sec:iota_A}
Let $\mathcal G$ be the coalescence of the  networks $\mathcal A$ and $\mathcal B$. Recall that the equilibrium points are given by the equation,
\[\mathscr E_{\mathbf y}= (\ker d)^\bot\cap  f^{-1}(\ker d),\]
for $\mathcal G$ and by
\[\mathscr E^{\mathcal A}_{\mathbf y}= (\ker d_{\mathcal A})^\bot\cap  f_{\mathcal A}^{-1}(\ker d_{\mathcal A}),\]
for ${\mathcal A}$, where the subindex ${\mathcal A}$ denotes the respective maps for the network ${\mathcal A}$  (similarly for ${\mathcal B}$). Now note that any generating cycle of $\ker d$ is entirely contained in ${\mathcal A}$ or in ${\mathcal B}$ as the only edges from ${\mathcal A}$ to ${\mathcal B}$ are the ones with $n'$ as a node. Thus,
\[\ker d = \iota_{\mathcal A}(\ker d_{\mathcal A})\oplus \iota_{\mathcal B}(\ker d_{\mathcal B}),\]
where, $\iota_{\mathcal A} (\mathbf y_{\mathcal A})=(\mathbf y_{\mathcal A}, \mathbf 0_{\mathcal B}) $ (resp. for $\iota_{\mathcal B}$). Moreover, it is clear that,
\[(\ker d)^\bot =\iota_{\mathcal A}((\ker d_{\mathcal A})^\bot)\oplus \iota_{\mathcal B}((\ker d_{\mathcal B})^\bot).\]
Now, as $f$ is defined in each component by one variable functions $f_{ij}$, we have $f(\mathbf y)=(f_{\mathcal A}(\mathbf y_{\mathcal A}),f_{\mathcal B}(\mathbf y_{\mathcal B}))$, and thus,
\[ f^{-1}(\ker d)=\iota _{\mathcal A}(f_{\mathcal A}^{-1}(\ker d_{\mathcal A}))\oplus \iota_{\mathcal B}( f_{\mathcal B}^{-1}(\ker d_{\mathcal B})).\]
In conclusion we have,
\begin{align*}
\mathscr E_{\mathbf y}=&\left (\iota_{\mathcal A} ((\ker d_{\mathcal A})^\bot)\oplus \iota_{\mathcal B} ((\ker d_{\mathcal B})^\bot )\right )\cap \left (\iota _{\mathcal A}( f_{\mathcal A}^{-1}(\ker d_{\mathcal A}))\oplus \iota_{\mathcal B}( f_{\mathcal B}^{-1}(\ker d_{\mathcal B}))\right )\\ 
=&\iota_{\mathcal A} \left ((\ker d_{\mathcal A})^\bot\cap f_{\mathcal A}^{-1}(\ker d_{\mathcal A})\right )\oplus \iota_{\mathcal B} \left ((\ker d_{\mathcal B})^\bot\cap  f_{\mathcal B}^{-1}(\ker d_{\mathcal B})\right  )\\
=& \iota_{\mathcal A}(\mathscr E_{\mathbf y} ^{\mathcal A})\oplus \iota_{\mathcal B}(\mathscr E_{\mathbf y} ^{\mathcal B}).
\end{align*}

\section{Phase portrait of a tree in low dimension}
\label{sec:pahse_portrait}
We study the qualitative dynamics of \cref{principal} in a line graph with three nodes where $f_{ij}=\hat f$ for all edges and all roots of $\hat f$ are simple. To simplify our study we also assume that the distance between consecutive roots of $\hat f$ are similar. That is, let $\alpha_0=0$, and $\{\alpha_i\}_{i= 1}^k$ be the positive roots of $\hat f$ in increasing order, where $k$ may be infinity. Then, we assume that,
\begin{equation}
    \max_{0\leq i<k} (\alpha_{i+1}-\alpha_i)\leq 2\min_{0\leq i<k} (\alpha_{i+1}-\alpha_i).
    \label{supp:eq:distance_alphas}
\end{equation}
If we let 1 be the middle node and 2, 3 the leaf nodes, the system \cref{principal} in  edge coordinates $y_1=y_{\{1,2\}}=x_2-x_1$ and $y_2=y_{\{1,3\}}=x_3-x_1$ is given by,
\begin{align}
\begin{split}
\dot y_1 &=  2\hat f(y_1)+\hat f(y_2),\\ 
\dot y_2 &=   \hat f(y_1)+2\hat f(y_2),
\end{split}
\label{eq:line_graph}
\end{align}
or compactly $\dot {\mathbf y}=g(\mathbf y)$. \Cref{eq:line_graph} is very symmetric, concretely if we consider a map that swaps the first and second coordinates, i.e. $p_{1,2}(y_1,y_2)=(y_2,y_1)$, we have, 
\[g\circ p_{1,2}=p_{1,2}\circ g,\hspace{1cm}\textrm{and}\hspace{1cm} g\circ -\textrm{id}=-\textrm{id}\circ g.\]
Thus the flow in $\mathbf y$ determines the flow in $\pm \mathbf y$ and $\pm (y_2,y_1)^\top$. In particular, we may restrict the study of the ODE in the region $y_2\geq |y_1|$. This situation immediately generalises when studying any star graph where we will have more possible permutations of variables. This is very reminiscent of the work done in  \cite{supp_MartinGolubitsky2006Ndon} as the symmetry of the star graph restricts the dynamics of the system.

Reparametrising time we can assume $\hat f'(0)>0$. Then, by  \cref{prop:tree_stability}, the equilibrium points of \cref{eq:line_graph} are of the form $(\pm\alpha_i, \pm\alpha _j)^\top$ and their stability is given by $\hat f'(\pm\alpha _i)$ and $\hat f'(\pm\alpha_j)$. Note that,
\[\hat f'(\pm\alpha _i)=\begin{cases}
>0		& $ if $ i $ is even$,\\
<0		& $ if $ i $ is odd$,\\
  \end{cases}\]
as all roots of $\hat f$ are simple so the signs alternate. We first compute,
\begin{align*}
\left.(\dot y_1+\dot y_2)\right|_{y_1+y_2=0} &=  \left.3(\hat f(y_1)+\hat f(y_2))\right|_{y_1+y_2=0} = 3(\hat f(y_1)+\hat f(-y_1))=0,\\ 
\left.(\dot y_1-\dot y_2)\right|_{y_1-y_2=0} &=  \left.(\hat f(y_1)-\hat f(y_2))\right|_{y_1-y_2=0} = \hat f(y_1)-\hat f(-y_1)=0,
\end{align*} 
so the lines $y_1=\pm y_2$ are invariant in our system.
\begin{figure}[htbp]
    \centering
        {\includegraphics[clip, trim=0cm 0.85cm 0cm 0.52cm, width=0.99\textwidth]{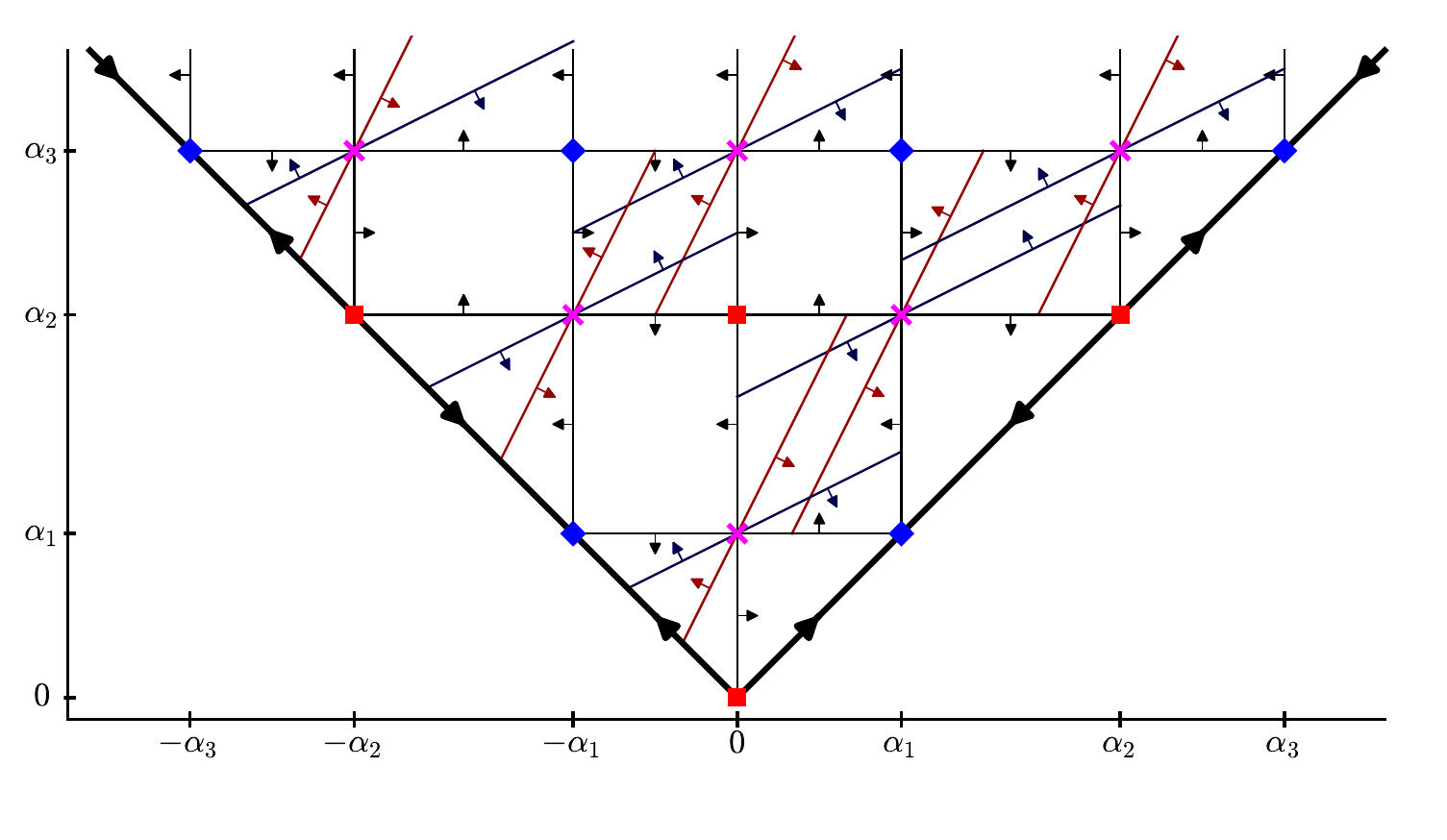}}
    \caption{Field \cref{eq:line_graph} in the invariant region $y_2\geq |y_1|$. \emph{Red squares}, \emph{blue diamonds} and \emph{magenta exes} represent respectively  source, sink and saddle  equilibrium points. The arrows represent the sense in which the field crosses the different color coded segments.}
    \label{fig:arrows2d}
\end{figure}
Moreover for all $i$,
\begin{align*}
\left.\dot y_1\right|_{y_1=\pm\alpha_i} &=  \left.(2\hat f(y_1)+\hat f(y_2))\right|_{y_1=\pm\alpha_i} = \hat f(y_2),\\ 
\left.\dot y_2\right|_{y_2=\pm\alpha_i} &=  \left.(\hat f(y_1)+2\hat f(y_2))\right|_{y_2=\pm\alpha_i} = \hat f(y_1),
\end{align*} 
so the horizontal and vertical segments between different equilibria are transversal to the field, i.e. orbits always cross them in the same sense, as depicted in black in \cref{fig:arrows2d}. Finally,
\begin{align*}
2\dot y_1 -\dot y_2 &= 3\hat f(y_1),\\ 
2\dot y_2- \dot y_1 &= 3\hat f(y_2),
\end{align*} 
and in particular the segments parallel to $2 y_1 = y_2$ (resp. $2 y_2 =y_1$), that go through the saddle equilibrium points, are transverse to the field for $y_1$ (resp. $y_2$) between roots, see in blue and red in  \cref{fig:arrows2d}. By \cref{supp:eq:distance_alphas}, these segments are contained in the regions delimited by the horizontal/vertical lines mentioned above.  

\begin{figure}[htbp]
    \centering
        {\includegraphics[clip, trim=0cm 0.85cm 0cm 0.44cm, width=0.99\textwidth]{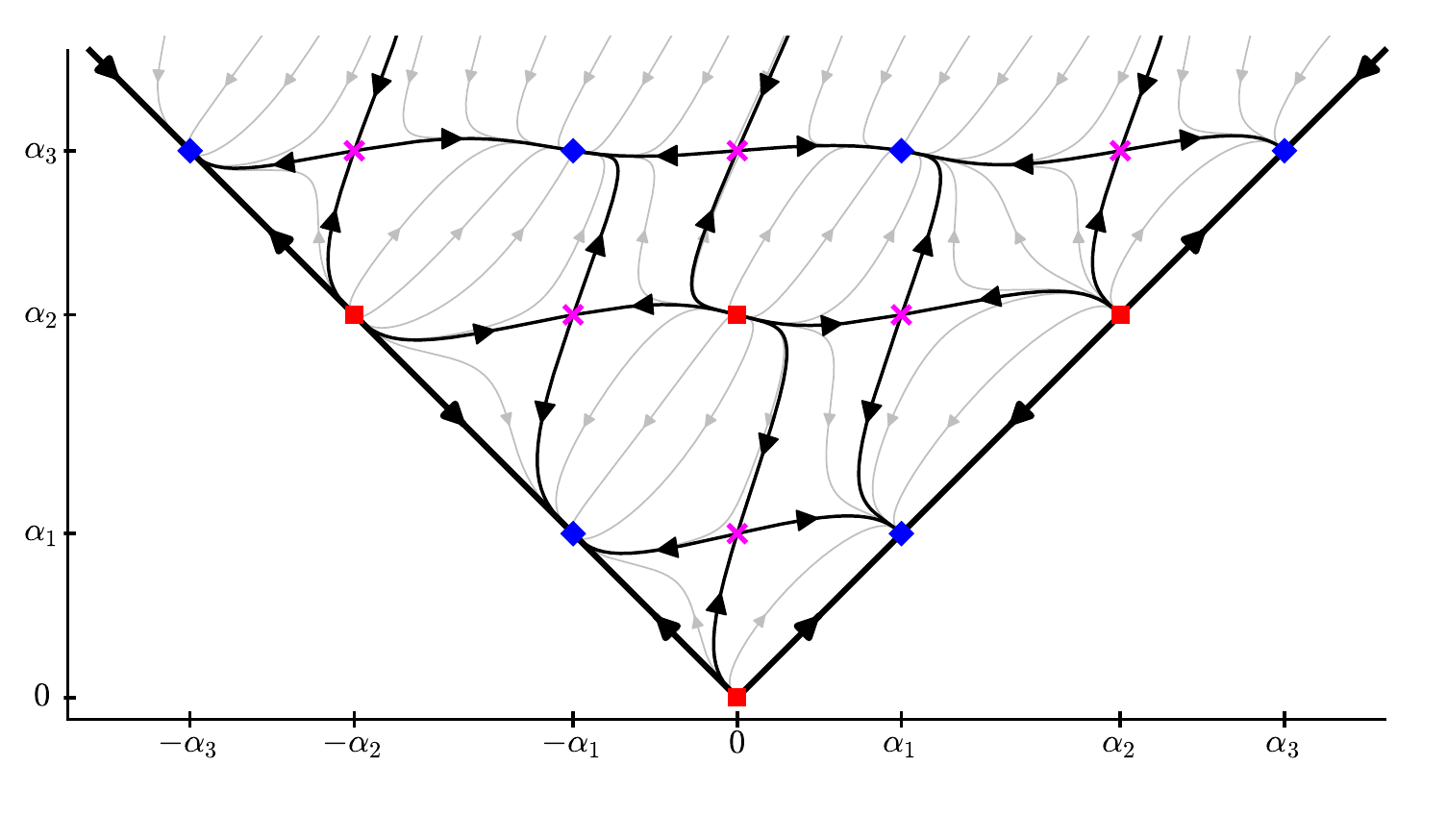}}
    \caption{Phase portrait of \cref{eq:line_graph} in the invariant region $y_2\geq |y_1|$. \emph{Red squares}, \emph{blue diamonds} and \emph{magenta exes} represent respectively  source, sink and saddles  equilibrium points. Lines in \emph{black} show the heteroclinic connections between saddles and other equilibrium points.}
    \label{fig:separatrix2d}
\end{figure}

It is not hard to see from \cref{fig:arrows2d}, that the  heteroclinic connections between the saddle and other equilibrium points are the ones depicted in \cref{fig:separatrix2d}.
 Moreover, by  \cref{prop:descomp_state_W} we know that the $\omega$ and $\alpha$-limit of any orbit is an equilibrium point or infinity, thus \cref{fig:separatrix2d} gives a complete description of the qualitative behaviour of the system.

\section{Equilibrium points for complete graphs}
\label{supp:sec:additive_open}
As in the main article we are working in the state space $\mathbb R^n/\langle\mathbf 1\rangle $ and we restrict our study to the equilibria  $\tilde {\mathscr E}_{\mathbf x}$ with $f_{i,j}=\hat f $ for all $\{i,j\}\in \mathcal E$. For each $\alpha \in \hat f ^{-1}(0)$, if $x_i^\star -x_j^\star\in \{0,\pm \alpha \}$ for all $i\not = j$, then $\mathbf x^\star \in \tilde {\mathscr E}_{\mathbf x}$. In this case, we may take the representative of $\mathbf x^\star$ such that $x_1^\star =0$. Then, for all $i$,
$x_i^\star =x_i^\star -x_1^\star \in \{0,\pm \alpha \}$.
Note that if $x_i=\alpha$ and $x_j=-\alpha $ for some $i,j$, then $x_i-x_j=\pm2\alpha $ which is a contradiction. Hence,  either $x^\star_i\in \{0,\alpha\}$ for all $i$ or $x^\star_i\in \{0,-\alpha\}$ for all $i$. In either case, choosing appropriate representative and ordering of the nodes,
 \begin{equation}
   \mathbf  x^\star =(0,\dots,0,\alpha,\dots,\alpha)^\top.
    \label{supp:(0...0.alpha...alpha)}
\end{equation}
 We show that generically these are all elements of $\tilde {\mathscr E}_{\mathbf x}$. Let $\mathbf x^\star \in \tilde {\mathscr E}_{\mathbf x} $, then $x_i^\star -x_j^\star \in \hat f ^{-1}(0)$  for all $i\not =j$. So if $i,j,k$ are pairwise different, and we denote $x_i^\star -x_j^\star=\alpha_1$, $x_j^\star -x_k^\star=\alpha_2$ and $x_k^\star -x_i^\star= \alpha _3$ we have $\alpha_1,\alpha_2,\alpha_3\in \hat f ^{-1}(0)$ and 
\begin{equation}
    \alpha_1+\alpha_2+\alpha_3=(x_i^\star -x_j^\star)+(x_j^\star -x_k^\star)+(x_k^\star -x_i^\star)=0.
    \label{eq:alpha1}
\end{equation}
 Now assume\footnote{This condition is generically satisfied as its negation implies a specific non-trivial linear relation between some roots of $\hat f$. } that $\hat f ^{-1}(0)\cap (0,\infty)$ is ``additive open'', i.e. $a+b\not =c$ for all $a,b,c\in \hat f ^{-1}(0)\cap (0,\infty)$. Note that in \cref{eq:alpha1} we may assume without loss of generality that $\alpha_3\leq 0$. Moreover, \cref{eq:alpha1} can be rewritten as $\alpha_1+\alpha_2=-\alpha_3$ where $\alpha_1, \alpha_2,-\alpha_3\geq 0$ and  as $\hat f$ is odd, $-\alpha_3$ is also a root of it. By the additive open condition, this can only be satisfied if $\alpha_l=0$ for some $l\in \{1,2,3\}$. Then, clearly the terms in the left hand of \cref{eq:alpha1} are  $0,\alpha,-\alpha$ for some root $\alpha$ of $\hat f$ which in principle could change between different triplets of nodes.
 Now fix $\alpha$ and $i, j$ two nodes such that $x_i^\star -x_j^\star=\pm\alpha$. Letting  $k$ run over all other nodes we conclude that
 $x_{i}^\star -x_{k}^\star \in \{0,\pm \alpha \}$ and $x_{j}^\star -x_{k}^\star \in \{0,\pm \alpha \}$. Repeating this process for other choices of $i,j$ we find that $x_{l_1}^\star -x_{l_2}^\star \in \{0,\pm \alpha \}$  for all $l_1\not =l_2$. Thus, by the argument made at the beginning of this section, $\mathbf x^\star$ can be represented by \cref{supp:(0...0.alpha...alpha)}.
 
  \section{Another deduction of the instability criterion in \texorpdfstring{$\mathcal K_n$}{Kn}}
 \label{supp:sec:instability_kn}
Let $\mathbf x^\star$ be as in \cref{supp:(0...0.alpha...alpha)}, denote by  $n_0$ the number of $0$'s of $\mathbf x^\star$ and by $n_\alpha=n-n_0$ the number of $\alpha$'s. Then, 
\[ L^+=\hat f'(0) \begin{pmatrix}n_0P_{n_0} & \mathbf 0_{n_0,n_\alpha}\\\mathbf 0_{n_\alpha,n_0} &n_\alpha P_{n_\alpha} \end{pmatrix} ,\hspace{18pt} L^- =|\hat f'(\alpha)|\begin{pmatrix}n_\alpha I_{n_0} & -\mathbf{1}_{n_0}\mathbf{1}_{n_\alpha}^\top\\-\mathbf{1}_{n_\alpha}\mathbf{1}_{n_0}^\top&n_0I_{n_\alpha} \end{pmatrix},\]
 where $P_{k}= I_k-\mathbf{1}_k\mathbf{1}_k^\top/k$, i.e. the orthogonal projection onto $\langle \mathbf 1 _k\rangle ^\bot$. 
Note that $L^+$ is the Laplacian of two disconnected complete graphs which we denote by $\mathcal K^+_0$ and $\mathcal K^+_\alpha$. Also note that $L^-$ is the Laplacian of a complete bipartite graph $\mathcal K_{n_0,n_\alpha}$. The effective resistances between nodes in these graphs are well known and can be deduced from \cite[Corollary C2]{supp_KleinD.1993Rd}. In our case we have,
\begin{equation*}
r_{ij}^+=\begin{cases}
\frac{2}{\hat f'(0)n_0}	& $ if $i,j\in \mathcal K^+_0,\\
\frac{2}{\hat f'(0)n_\alpha} & $ if $i,j\in \mathcal K^+_\alpha,\\
\infty & $else$,
  \end{cases}\hspace{1cm} \textrm{and}\hspace{1cm} r_{ij}^-=\begin{cases}
\frac{2}{|\hat f'(\alpha)|n_\alpha}	& $ if $i,j\in \mathcal K^+_0,\\
\frac{2}{|\hat f'(\alpha)|n_0} & $ if $i,j\in \mathcal K^+_\alpha,\\
\not =\infty & $else$. 
\end{cases}
\end{equation*}
By \cref{thm:unstability_r^-_ij r^+_ij}, we have instability if $r_{ij}^->r_{ij}^+$ for some edge $\{i,j\}$, which in or context reduces to, 
\[\frac{n_0}{n_\alpha}< \frac{\hat f'(0)}{|\hat f'(\alpha)|}\hspace{0.8cm}\textrm{or}\hspace{0.8cm}\frac{n_0}{n_\alpha}> \frac{|\hat f'(\alpha)|}{\hat f'(0)},\] 
or equivalently,
\begin{equation}
\frac{n_0}{n_\alpha}\not \in \left [  \frac{\hat f'(0)}{|\hat f'(\alpha)|},\frac{|\hat f'(\alpha)|}{\hat f'(0)}\right ].
\label{supp:eq:in[]}
\end{equation} 
Now we need the following elementary relation between inequalities. Given positive numbers $a,b,c,d>0$, then 
\[\frac{a}{b}>\frac{c}{d} \Leftrightarrow \frac{b}{a}<\frac{d}{c} \Leftrightarrow \frac{b}{a}+1<\frac{d}{c}+1 \Leftrightarrow \frac{a+b}{a}<\frac{c+d}{c} \Leftrightarrow \frac{a}{a+b}>\frac{c}{c+d},\]
and similarly,
\[\frac{a}{b}<\frac{c}{d} \Leftrightarrow  \frac{a}{a+b}<\frac{c}{c+d}.\]
Thus, condition \cref{supp:eq:in[]} is equivalent to,
\[\frac{n_0}{n_0+n_\alpha}=\frac{n_0}{n}\not \in \left [  \frac{\hat f'(0)}{\hat f'(0)+|\hat f'(\alpha)|},\frac{|\hat f'(\alpha)|}{\hat f'(0)+|\hat f'(\alpha)|}\right ],\] 
which is precisely the one we got in \cref{prop:kn_stability} by computing the eigenvalues, and hence it is also tight.

\bibliographystyle{siamplain}
\bibliography{references_supp}